\documentclass[leqno, a4paper]{amsart}
\usepackage{amssymb,latexsym}

\author[W.\thinspace{}Kim]{Wansu Kim}
\address{Wansu Kim\\%
Department of Mathematics\\%
South Kensington Campus\\%
Imperial College London\\%
London, SW7 2AZ\\%
United Kingdom}
\email{w.kim@imperial.ac.uk}

\usepackage[active]{srcltx}
\usepackage{amsmath}
\usepackage{amsfonts}

\usepackage{amsthm}
\usepackage{amscd}
\usepackage{mathrsfs}
\usepackage{fancyhdr}
\usepackage{graphicx}
\usepackage{psfrag}
\usepackage{calc}
\usepackage{url}

\usepackage{pb-diagram,pb-xy}
\usepackage{stmaryrd}
\usepackage{times}
\usepackage{verbatim}
\usepackage[cmtip,arrow,matrix,curve,tips,frame]{xy}
\usepackage{hyperref}
\xyoption{matrix}

\usepackage{ifthen}
\usepackage{calc}
\usepackage{epsfig}
\usepackage{color}
\usepackage{amsxtra}
\usepackage{amstext}
\usepackage{amssymb}
\usepackage{stmaryrd}
\usepackage{mathrsfs}
\usepackage{pifont}
\usepackage{charter}
\usepackage{avant}
\usepackage{courier}
\usepackage[T1]{fontenc}
\usepackage{textcomp}

\usepackage{egastyle}
\numberwithin{equation}{subsection}

\theoremstyle{plain}
\newtheorem{thm}[subsection]{Theorem}
\newtheorem*{thm*}{Theorem}

\newtheorem{thmsub}[equation]{Theorem}

\newtheorem*{exthm*}{Expected Theorem}

\newtheorem{lem}[subsection]{Lemma}
\newtheorem{lemsub}[equation]{Lemma}
\newtheorem*{lem*}{Lemma}

\newtheorem{prop}[subsection]{Proposition}

\newtheorem{propsub}[equation]{Proposition}
\newtheorem*{prop*}{Proposition}

\newtheorem{cor}[subsection]{Corollary}
\newtheorem{corsub}[equation]{Corollary}
\newtheorem*{cor*}{Corollary}

\newtheorem{claimsub}[equation]{Claim}
\newtheorem*{claim*}{Claim}
\newtheorem*{nclma}{Claim A}

\newtheorem{conj}[subsection]{Conjecture}
\newtheorem*{conj*}{Conjecture}

\theoremstyle{definition}
\newtheorem{defn}[subsection]{Definition}
\newtheorem*{defn*}{Definition}
\newtheorem{defnsub}[equation]{Definition}

\newtheorem{exasub}[equation]{Example}
\newtheorem*{exa*}{Example}

\theoremstyle{remark}
\newtheorem{rmk}[subsection]{Remark}
\newtheorem{rmksub}[equation]{Remark}
\newtheorem*{rmk*}{Remark}

\numberwithin{figure}{subsection}
\numberwithin{table}{subsection}

\newcounter{listnum}

\DeclareMathOperator{\coker}{coker}

\DeclareMathOperator{\rank}{rank}

\DeclareMathOperator{\id}{id}

\newcommand{\e}{\mathbf{e}}

\DeclareFontEncoding{OT2}{}{} % to enable usage of cyrillic fonts

\newcommand{\nf}[1]{\underline{#1}}
\newcommand{\ol}[1]{\overline{#1}}
\newcommand{\wt}[1]{\widetilde{#1}}

\newcommand{\tim}{\!\cdot\!}
\newcommand{\geqs}{\geqslant}
\newcommand{\leqs}{\leqslant}

\newcommand{\et}{\text{\rm\'et}}
\DeclareMathOperator{\cyc}{cyc}

\newcommand{\wh}[1]{\widehat{#1}}

\newcommand{\M}{\mathcal{M}}

\newcommand{\llpar}{(\!(}
\newcommand{\rrpar}{)\!)}

\newcommand{\rep}{\mathfrak{R}}

\DeclareMathOperator{\alg}{alg}
\DeclareMathOperator{\sep}{sep}

\newcommand{\starr}{^\times}

\DeclareMathOperator{\Frac}{Frac}

\newcommand{\ra}{\rightarrow}
\newcommand{\xra}[1]{\xrightarrow{#1}}

\newcommand{\hra}{\hookrightarrow}

\newcommand{\thra}{\twoheadrightarrow}

\newcommand{\lra}{\longrightarrow}

\newcommand{\self}[1]{{#1}\ra {#1}}

\newcommand{\invlim}{\mathop{\varprojlim}\limits}

\newcommand{\Gm}{\mathbb{G}_m}

\newcommand{\set}[1]{\{#1\}}
\newcommand{\iv}{^{-1}}
\newcommand{\ivtd}[1]{\frac{1}{#1}}

\newcommand{\Eps}{\mathcal{E}}
\newcommand{\vphi}{\varphi}

\newcommand{\Sig}{\mathfrak{S}}

\newcommand{\Q}{\mathbb{Q}}

\newcommand{\cD}{\mathcal{D}}

\newcommand{\PP}{\mathcal{P}}

\newcommand{\kbar}{\bar{k}}

\newcommand{\Z}{\mathbb{Z}}

\newcommand{\F}{\mathbb{F}}

\newcommand{\Qp}{\Q_p}

\newcommand{\Kbar}{\overline{K}}

\newcommand{\Zp}{\Z_p}

\newcommand{\Fp}{\F_p}

\newcommand{\fo}{\mathscr{O}}

\DeclareMathOperator{\Gal}{Gal}
\newcommand{\gal}{\boldsymbol{\mathcal{G}}}
\newcommand{\GK}{\gal_K}
\newcommand{\GKinfty}{\gal_{K_\infty}}

\newcommand{\fa}{\mathfrak{a}}

\newcommand{\m}{\mathfrak{m}}

\newcommand{\gM}{\mathfrak{M}}
\newcommand{\gN}{\mathfrak{N}}

\DeclareMathOperator{\ur}{ur}

\DeclareMathOperator{\Hom}{Hom}

\DeclareMathOperator{\Tor}{Tor}

\DeclareMathOperator{\Fil}{Fil}
\DeclareMathOperator{\gr}{gr}

\DeclareMathOperator{\Mod}{Mod}

\DeclareMathOperator{\cris}{cris}
\newcommand{\Acris}{A_{\cris}}

\DeclareMathOperator{\st}{st}

\DeclareMathOperator{\dR}{dR}

\DeclareMathOperator{\free}{free}

\newcommand{\phimod}[1]{\nf{\Mod}_{#1}(\varphi)}
\newcommand{\phimodh}[1]{\nf{\Mod}_{#1}(\varphi)^{\leqs h}}
\newcommand{\phimodBT}[1]{\nf{\Mod}_{#1}(\varphi)^{\leqs 1}}

\DeclareMathOperator{\Rep}{Rep}

\newcommand{\freeprep}{\Rep_{\Zp}^{\free}}

\newcommand{\Repcrh}{\Rep_{\cris,\Qp}^{[0,h]}(\GK)}
\newcommand{\Repsth}{\Rep_{\st,\Qp}^{[0,h]}(\GK)}

\newcommand{\Replatcrh}{\Rep_{\cris,\Zp}^{[0,h]}(\GK)}
\newcommand{\Replatsth}{\Rep_{\st,\Zp}^{[0,h]}(\GK)}

\newcommand{\DD}{\mathbb{D}}

\title[Classification of $p$-divisible groups]{The classification of $p$-divisible groups over $2$-adic discrete valuation rings}

\begin{document} 

\begin{abstract}
Let $\fo_K$ be a $2$-adic discrete valuation ring with perfect residue field $k$. We classify $p$-divisible groups and $p$-power order finite flat group schemes over $\fo_K$ in terms of certain Frobenius modules over $\Sig:=W(k)[[u]]$.  We also show the compatibility with crystalline Dieudonn\'e theory and associated Galois representations. Our approach differs from Lau's generalization of display theory, who independently obtained our result using display theory.
\end{abstract}

\maketitle
\tableofcontents
\section{Introduction}\label{sec:intro}
Let $k$ be a perfect field of characteristic $p$, $W(k)$ its ring of Witt vectors, and $K_0:=W(k)[\ivtd p]$. Let $K$ be a finite totally ramified extension of $K_0$ and let us fix its algebraic closure $\Kbar$. We fix a uniformizer $\pi\in K$ and an Eisenstein polynomial $\PP(u)\in W(k)[u]$ such that $\PP(\pi)=0$ and $\PP(0)=p$. Choose $\pi^{(n)}\in\Kbar$ so that $(\pi^{(n+1)})^p = \pi^{(n)}$ and $\pi^{(0)}=\pi$. Put $K_\infty:= \bigcup_n K(\pi^{(n)})$, $\gal_K:=\Gal(\Kbar/K)$, and $\gal_{K_\infty}:=\Gal(\Kbar/K_\infty)$. 

Set $\Sig:=W(k)[[u]]$ and we extend the Witt vector Frobenius map to a ring endomorphism $\vphi:\Sig\ra\Sig$ by $\vphi(u)=u^p$. Breuil has conjectured that $p$-divisible groups and finite flat group schemes over $\fo_K$ are classified by a certain Frobenius-module over $\Sig$; namely, $\phimodBT\Sig$ and $(\Mod/\Sig)^{\leqs1}$, respectively. (The precise definitions will be given later in  Definitions~\ref{def:KisinLatt} and \ref{def:KisinTor}.)  Also note that there exist contravariant functors $  T^*_\Sig$ from  $\phimodBT\Sig$ and $(\Mod/\Sig)^{\leqs1}$ to the category of finitely generated $\Zp$-modules with continuous $\GKinfty$-action (\ref{eqn:TSigLattice}, \ref{eqn:TSigTors}). The main result of this paper is the proof of the following conjecture of Breuil.
\begin{conj}[Breuil]\label{conj:breuil}
There exists an exact anti-equivalence of categories $G\leftrightsquigarrow\gM_G$ between the category of $p$-divisible groups over $\fo_K$ and $\phimodBT\Sig$ such that the contravariant Dieudonn\'e crystal $\DD^*(G)$ can be recovered from $\gM_G$ and  there exists a natural $\GKinfty$-equivariant isomorphism $T_p(G)\cong T^*_\Sig(\gM_G)$. 

There exists an exact anti-equivalence of categories $H\leftrightsquigarrow\gM_H$ between the category of $p$-power order finite flat group schemes over $\fo_K$ and $(\Mod/\Sig)^{\leqs1}$ such that the contravariant Dieudonn\'e crystal $\DD^*(H)$ can be recovered from $\gM_H$ and  there exists a natural $\GKinfty$-equivariant isomorphism $H(\Kbar)\cong T^*_\Sig(\gM_H)$. 
%
%Consider a sequence $G^\bullet:\ 0\ra G' \ra G \ra G'' \ra 0$ where $G$, $G'$, and $G''$ are either $p$-power order finite flat group schemes or $p$-divisible groups over $\fo_K$ and let $\gM^\bullet$ denote the corresponding sequence of  $(\vphi,\Sig)$-modules. Then $G^\bullet$ is short exact if and only if $\gM^\bullet$ is short exact. (The statement also holds when $G'$ is a finite flat group scheme and $G$ and $G''$ are $p$-divisible groups.)
\end{conj}
The definition of these functors make use of the choice of uniformizer $\pi\in\fo_K$ (since the definitions of the target categories $\phimodBT\Sig$ and $(\Mod/\Sig)^{\leqs1}$ do).\footnote{Note that for fixed $\pi=\pi^{(0)}$ any two choices of $\GKinfty\subset \GK$ are conjugate to each other.}  Note that the $\GKinfty$-restriction  induces a fully faithful functor from the category of $\GK$-representations coming from either a $p$-divisible group or a finite flat group scheme over $\fo_K$ to the category of $\GKinfty$-representations (\emph{cf.} \cite[Corollary~2.1.14]{kisin:fcrys} and Corollary~\ref{cor:FF} in this paper), and  when $G$ is a $p$-divisible group we can give a description of the $\GK$-action on $T^*_\Sig(\gM_G)$ which recovers the natural $\GK$-action on $T_p(G)$. (\emph{cf.} \eqref{eqn:TstGK}, Proposition~\ref{prop:main2}, and the proof of Proposition~\ref{prop:key}.) 

Although the compatibility of the functors with the associated Galois representations was not explicitly stated in the original conjecture \cite[\S2]{breuil:GpSchNormField}, it  is important for the applications in number theory; for example, this is needed for the construction of moduli of finite flat group schemes \cite{Kisin:ModuliGpSch},  and Kisin's strategy of constructing a good integral  model of an abelian-type Shimura variety \cite{Kisin:IntModelAbType}.\footnote{Note that Vasiu  \cite{Vasiu:GoodRedn2} has claimed to have proved the main result of \cite{Kisin:IntModelAbType} including  $2$-adic case without using Conjecture~\ref{conj:breuil}.}  In the latter case we need this extra compatibility because Kisin's argument requires that the assignment $T_p(G)\rightsquigarrow \gM_G$ for a $p$-divisible group $G$ over $\fo_K$ should extend to a $\otimes$-functor on the $\otimes$-category generated by the $T_p(G)$'s in the category of $\Zp[\GK]$-modules (hence, we escape the realm of Barsotti-Tate representations), and this could be checked via Galois compatibility.

Now we comment on the previous results on the conjecture. Breuil \cite{breuil:GpSchNormField, Breuil:GrPDivGrFiniModFil} proved the conjecture for finite flat group schemes killed by $p$ when $p>2$. Kisin  \cite[\S{A} and Corollary~2.2.6]{kisin:fcrys} proved the full conjecture when $p>2$, and proved the ``isogeny version'' of the conjecture when $p=2$. %(See Theorem\ref{thm:Kisin}(\ref{thm:Kisin:equiv}) for the precise statement.) 
%; i.e., Kisin constructed an anti-equivalence between the isogeny category of $p$-divisible groups over $\fo_K$ and the isogeny category $\phimodBT\Sig[\ivtd p]$ which is compatible with the associated $\GKinfty$-representations and the Dieudonn\'e isocrystals (or filtered $\vphi$-module, if you wish). 
Later Kisin \cite[\S1]{Kisin:2adicBT} proved the conjecture for the Cartier duals of connected $p$-divisible groups without restriction on $p$, using Zink's theory of windows and displays. In this paper, we prove the conjecture without connectedness assumptions when $p=2$.

Independently, Eike~Lau \cite{Lau:2010fk, Lau:GalRep} has proven the conjecture by a different approach. Indeed, Lau \cite{Lau:2010fk} extended the Zink's theory of windows and displays for arbitrary $p$-divisible groups with no restrictions on $p$ (over more general bases than discrete valuation rings), and obtained the  compatibility with Galois representations in \cite{Lau:GalRep} via ``displayed equation''. Shortly after this paper was completed, Tong~Liu \cite{Liu:Classif} independently proved Theorem~\ref{thm:classif} of this paper using his theory of $(\vphi,\hat G)$-modules. Our approach is different from aforementioned results, and more related to Kisin's original approach in   \cite{kisin:fcrys}. 

The $p=2$ case of Conjecture~\ref{conj:breuil} has been the main technical obstacle to constructing a good $2$-adic integral model of a Shimura variety of abelian type following Kisin's strategy \cite{Kisin:IntModelAbType}. Also Madapusi Pera's work on toroidal compactifications of these integral models in the Hodge-type case \cite{MadapusiPera:Thesis} assumes $p>2$ for the similar reason. Though we have not succeeded yet, it seems that the result of this paper (which proves Conjecture~\ref{conj:breuil}) can be used to extend the main results of \cite{Kisin:IntModelAbType} (and quite possibly, \cite{MadapusiPera:Thesis})  to cover some $2$-adic case. (See Remark~\ref{rmk:IntModel}  for more details.) Note that the existence of such good integral models plays a fundamental role in Langlands' Jugendtraum \cite{Langlands:Jugendtraum} and Kudla's program outlined in \cite{Kudla:EisensteinOutline}.

Now, let us explain some of the trouble one encounters when proving Conjecture~\ref{conj:breuil} when $p=2$ (with no connectedness assumptions). We consider the following the $\Zp$-lattices in $V_p(G):=T_p(G)[\ivtd p]$ for a $p$-divisible group $G$ over $\fo_K$.
\begin{itemize}
\item $T:= T_p(G)$, the $p$-adic Tate module.
\item $T_\gM:=T^*_\Sig(\gM_G)$, where $T^*_\Sig(\gM_G)$ is defined in \eqref{eqn:TSigLattice}.
\item $T_\DD$, the $\Zp$-lattice associated to the Dieudonn\'e crystal via  $T^*_{\st}$ \eqref{eqn:Tst}. \emph{cf.} \S\ref{subsec:Dieudonne}.
\end{itemize}
The desired compatibility follows from the equality $T=T_\gM$.

When $p>2$ or when the Cartier dual $G^\vee$ is connected,  we can prove the equality $T=T_\gM$ by showing $T=T_\DD$ and $T_\gM=T_\DD$. Indeed, the equality $T_\gM=T_\DD$ follows from a standard linear algebra computation together with the compatibility with crystalline Dieudonn\'e theory, if $p>2$ or if the Cartier dual $G^\vee$ is connected. The equality $T=T_\DD$ follows from Faltings's integral comparison isomorphism when $p>2$ \cite[Theorem~7]{Faltings:IntegralCrysCohoVeryRamBase} and Zink's theory of windows and display when $G^\vee$ is connected \cite[Proposition~1.1.10]{Kisin:2adicBT}. 
The trouble is that when $p=2$ the lattice $T_\gM$ is \emph{strictly} contained in $T_\DD$ if $G^\vee$ is not connected (Proposition~\ref{prop:main2}). Indeed, we will show  that $T=T_\gM$ as sublattices in $T_\DD$.

Let us explain our approach to prove Conjecture~\ref{conj:breuil}. We use Kisin theory \cite{kisin:fcrys} to define a functor $G\rightsquigarrow\gM_G$ from the category of $p$-divisible groups to the category of $(\vphi,\Sig)$-modules in such a way that the equality $T=T_\gM$ is automatic by construction. We show that this functor is an anti-equivalence, and then we separately show  the compatibility with Dieudonn\'e crystals. Since the result up to isogeny was already obtained by Kisin \cite{kisin:fcrys}, it is not surprising that the integral results can be obtained by studying Galois-stable lattices in $V_p(G)$.

The author expects that this idea can be generalized to show the following statement about $\GK$-stable $\Zp$-lattices in semi-stable representations with small Hodge-Tate weights.
\begin{conj}
Let $D:=(D,\vphi_D,N,\Fil^\bullet D_K)$ be a weakly admissible filtered $(\vphi,N)$-module such that $\gr^wD_K = 0$ for all $w\notin [0,p-1]$. Consider $\gM\in\phimod\Sig^{\leqs p-1}$ equipped with a $\GKinfty$-equivariant embedding $T^*_\Sig(\gM)\hra V^*_{\st}(D)$. (See Definition~\ref{def:KisinLatt} for the definition of $\phimod\Sig^{\leqs p-1}$.) Then the image of the embedding is $\GK$-stable if and only if $S\otimes_{\vphi,\Sig}\gM$ is a strongly divisible $S$-lattice in $S\otimes_{W(k)}D$ in the sense of \cite[Definition~2.2.1]{Breuil:IntegralPAdicHodgeThy}. Furthermore, if $D$ does not admit a non-zero weakly admissible quotient pure of slope $p-1$, then we have $T^*_\Sig(\gM) = T^*_{\st}(S\otimes_{\vphi,\Sig}\gM)$ as $\Zp$-lattices in $V^*_{\st}(D)$.
\end{conj}
This conjecture is proved  by Tong Liu \cite{Liu:StronglyDivLattice} when $\gr^wD_K = 0$ for all $w\notin [0,p-2]$, so the only unknown case is when $\gr^0 D_K\ne 0 $ and $\gr^{p-1}D_K\ne 0$.

We include sections to review Kisin theory (\S\ref{sec:KisinThy}) and various definitions and constructions involving Breuil modules (\S\ref{sec:DefsBreuil}). The main results are stated in \S\ref{sec:statements}, and proved in \S\ref{sec:proofs}. We allow $p$ to be arbitrary prime in order to illustrate the differences when $p=2$, though our proof is not essentially different from Kisin's \cite{kisin:fcrys, Kisin:2adicBT} when $p>2$.
\subsection*{Acknowledgement}
The author deeply thanks Brian Conrad, Mark Kisin, and Tong Liu for helpful comments. Also he especially thanks the anonymous referee who gave a very careful reading to this paper and provided numerous helpful suggestions, especially on the proof of Proposition~\ref{prop:key}.

\section{Review: Kisin theory}\label{sec:KisinThy}
Let $\Sig:=W(k)[[u]]$ where $u$ is a formal variable. Let $\fo_\Eps$ be the $p$-adic completion of $\Sig[\ivtd u]$, and set $\Eps:=\fo_\Eps[\ivtd p]$. Note that $\fo_\Eps$ is a complete discrete valuation ring with uniformiser $p$ and $\fo_\Eps/(p)\cong k\llpar u \rrpar$.  %We view $k\llpar u \rrpar$ as the norm field $ X_K(K_\infty)$  for the extension $K_\infty /K$ \cite[\S2]{wintenberger:NormFiels}.\footnote{Note that $X_K(K_\infty)$ is a local field of characteristic $p$ with residue field $k$ \cite[Proposition~2.2.4]{wintenberger:NormFiels}, and one can directly check that $\nf\pi$ is the uniformizer of $ X_K(K_\infty)$ by the definition of the valuation \cite[\S2.2.3.3]{wintenberger:NormFiels}.} 
We extend the Witt vectors Frobenius map to $\Sig$, $\fo_\Eps$, and $\Eps$ by sending $u$ to $u^p$, and denote them by $\vphi$. %(We write $\vphi_\Sig$ instead, if we need to specify that it is an endomoprhism on $\Sig$, for example.) 
%Note that $\vphi$ is finite and flat. 
We denote by $\vphi^*(\cdot)$ the scalar extension by $\vphi$. We fix an Eisenstein polynomial $\PP(u)\in W(k)[u]$ such that $\PP(\pi)=0$ and $\PP(0)=p$, and view it as an element of $\Sig$. Note that this normalization is different from the usual convention which takes $\PP(u)$ to be a monic polynomial -- although our discussions does not depend upon the choice of normalization of $\PP(u)$ in any significant way, our convention seems more natural in some places. (See \S\ref{par:DefEtPhiNilp}, for example.)

\begin{defn}\label{def:KisinLatt}
For a non-negative integer $h$, a  \emph{$\vphi$-module of height $\leqs h$} is a tuple $(\gM,\vphi_\gM)$, where $\gM$ is a finite free $\Sig$-module, and  $\vphi_\gM:\gM\ra\gM$ is a $\vphi$-semilinear endomorphism such that $\PP(u)^h$ annihilates $\coker(1\otimes\vphi_\gM:\vphi^*\gM \ra \gM)$. Let $\phimodh\Sig$ denote the category of $\vphi$-modules of height $\leqs h$ with the obvious notion of morphisms. When $h=1$, we set $\phimodBT\Sig:=\Mod_\Sig^{\leqs1}(\vphi)$.
\end{defn}

\subsection{Galois representations}\label{subsec:GKinftyRep} 
Put $\rep:=\invlim_{x\mapsto x^p} \fo_{\ol K }/(p)$. Note that  the sequence $\nf\pi:=\set{\pi^{(n)}}$ is an element of $\rep$, where $\pi^{(n)}\in K_\infty$ is the fixed $p^n$th root of the fixed uniformizer $\pi=\pi^{(0)}$. (See the beginning of \S\ref{sec:intro}.) It follows from \cite[Corollaires~3.2.3, 4.3.4]{wintenberger:NormFiels} that $\Frac\rep$ is isomorphic to the $\nf\pi$-adic completion of $k\llpar \nf\pi \rrpar^{\alg}$.\footnote{To compare the notation of \cite{wintenberger:NormFiels} with ours, $R(\ol K) = \Frac\rep$ and $ X_K(K_\infty) = k\llpar \nf\pi \rrpar$; to see the second equality, note that $ X_K(K_\infty)$ is a local field of characteristic $p$ with residue field $k$ \cite[Proposition~2.2.4]{wintenberger:NormFiels}, and one can directly check that $\nf\pi$ is the uniformizer of $ X_K(K_\infty)$ by the definition of the valuation \cite[\S2.2.3.3]{wintenberger:NormFiels}.} Also $\rep$ is the valuation ring for the $\nf\pi$-adic  valuation on $\Frac\rep$ (\emph{cf.} \cite[\S4.1.1]{wintenberger:NormFiels}).

We have a $W(k)$-embedding $\Sig\hra W(\rep)$ defined by sending $u$ to $[\nf\pi]$, where $[\cdot]$ is the Teichm\"uller lift. This embedding extends to $\fo_\Eps \hra W(\Frac\rep)$. Let $\wh\fo_{\Eps^{\ur}}$ be the topological closure of the integral closure of $\fo_{\Eps}$ in $W(\Frac\rep)$, and $\wh\Sig^{\ur}$ the topological closure of the integral closure of $\Sig$ in $\wh\fo_{\Eps^{\ur}}$. (We use the $p$-adic topology in both cases.) Clearly, $\wh\Sig^{\ur}$ and $\wh\fo_{\Eps^{\ur}}$ are subrings of $W(\Frac\rep)$ stable under $\vphi$ and $\GKinfty$-action (but \emph{not} under $\GK$-action), and $\wh\Sig^{\ur}$ is contained in $W(\rep)$.

Now we associate, to $\gM\in\phimodh\Sig$, a $\Zp$-module with continuous $\GKinfty$-action, as follows:
\begin{equation}\label{eqn:TSigLattice}
T^*_\Sig(\gM):=\Hom_{\Sig,\vphi}(\gM,\wh\Sig^{\ur}),
\end{equation}
where $\GKinfty$ acts through $\wh\Sig^{\ur}$. By Proposition~1.8.3 and \S{A} 1.2 of \cite{fontaine:grothfest}, it follows that  $T^*_\Sig(\gM)$ is $\Zp$-free with rank equal to the $\Sig$-rank of $\gM$, and  $T_\Sig^*$ is exact and commutes with direct sums and $\otimes$-products. 

For a topological group $\Gamma$ (e.g., $\GK$, $\GKinfty$), let $\freeprep(\Gamma)$ denote the category of finite free $\Zp$-modules with continuous $\Gamma$-action. Let $\Repcrh$ and $\Repsth$ respectively denote the categories of crystalline and semi-stable $\GK$-representations $V$ such that  $\gr^wD_{\dR}^*(V) = 0$ for $w\notin[0,h]$. Let $\Replatcrh$ (respectively, $\Replatsth$) denote the full subcategory of $\freeprep(\GK)$ consisting of objects $T$ such that $T[\ivtd p]\in \Repcrh$ (respectively, $T[\ivtd p]\in \Repsth$).
The following  non-trivial theorem of Kisin is the starting point of our argument:

\begin{thm}[Kisin]\label{thm:Kisin}\hfill
\begin{enumerate}
\item
\cite[Proposition~2.1.12]{kisin:fcrys}
The contravariant functor $T^*_\Sig:\phimodh\Sig\ra\freeprep(\GKinfty)$ is fully faithful.
\item
\label{thm:Kisin:main}
There exists a (contravariant) functor  $\gM^*:\Replatsth\ra\phimodh\Sig$ of $\otimes$-categories such that for any $T\in\Replatsth$ there exists a natural $\GKinfty$-isomorphism $T\cong T^*_\Sig(\gM^*(T))$.
 \item 
 \label{thm:Kisin:equiv}
 \cite[Theorem~2.2.7]{kisin:fcrys}
 Consider the following functor\footnote{The functor can be defined because $V:= V_p(G)$ is a crystalline representation with $\gr^w D^*_{\dR}(V)=0$ for $w\ne0,1$.}
 \begin{equation*}\label{eqn:gM}
\gM^*\circ T_p :
\left\{p\text{-divisible groups over }\fo_K\right\} \ra \phimodBT\Sig;\ G \rightsquigarrow \gM^*(T_p(G)).
\end{equation*}
This functor is an  anti-equivalence of categories if $p>2$, and induces an  anti-equivalence of categories between the corresponding isogeny categories if $p=2$.
%The functor $(\gM^*\circ T_p)[\ivtd p]:[G]\rightsquigarrow \gM^*(T_p(G))[\ivtd p]$ from the isogeny category of $p$-divisible groups over $\fo_K$ to the isogeny category $\phimodBT\Sig[\ivtd p]$ is an exact anti-equivalence of categories.\footnote{Clearly, the definition of $[G]\rightsquigarrow \gM^*(T_p(G))[\ivtd p]$ does not depend on the choice of $G$ in its isogeny class.}
\end{enumerate}
\end{thm}
\begin{proof}
We only explain how to read off \eqref{thm:Kisin:main} from \cite{kisin:fcrys}. Kisin constructed the functor
$\Repsth \ra \phimodh\Sig[\ivtd p]$ (see \cite[Corollary~1.3.15]{kisin:fcrys}), and showed  that if this functor sends $V\in\Repsth$ to the isogeny class containing $\gM'\in\phimodh\Sig$ then we have a natural $\GKinfty$-isomorphism $V\cong T^*_\Sig(\gM')[\ivtd p]$. (See \cite[Proposition~2.1.5]{kisin:fcrys}.) Finally, the proof of  \cite[Lemma~2.1.15]{kisin:fcrys} shows that for any $\GKinfty$-stable $\Zp$-lattice $T$ in  $V\cong T^*_\Sig(\gM')[\ivtd p]$ there exists a unique $\vphi$-stable $\Sig$-lattice $\gM\subset \gM'[\ivtd p]$ which is of height $\leqs h$, such that the $\Zp$-lattice $T^*(\gM)$ in $V$ is precisely $T$. Set $\gM^*(T):= \gM$.
\end{proof}
\begin{rmk}\label{rmk:PathologyLiu}
The functor $\gM^*$ is unfortunately \emph{not} exact although it is exact up to isogeny.\footnote{The author thanks Tong~Liu for pointing this out to him.} In fact, T.~Liu \cite[Example~2.5.6]{Liu:LattFiltPhiNMod}  provided an example of a short exact sequence $T^\bullet$ of Tate modules of $p$-divisible groups such that $\gM^*(T^\bullet)$ is not exact. In that example $T^\bullet$ does \emph{not} extend to a short exact sequence of $p$-divisible groups over $\fo_K$,\footnote{Although by Tate's theorem \cite[(4.2)]{Tate:pDivGps} an exact sequence $0\ra T_p(G')\ra T_p(G) \ra T_p(G'')\ra 0$ extends to a sequence $0\ra G' \ra G \ra G'' \ra 0$ of $p$-divisible groups,  this sequence does not have to be exact.} and we will prove later (\S\ref{subsec:exactness}) that $\gM^*$ preserves a short exact sequence of Tate modules whenever it extends to a short exact sequence of $p$-divisible groups. 
\end{rmk}
%\begin{rmk}\label{rmk:KisinBT}
%We sketch the construction of the quasi-inverse $\gM[\ivtd p]\rightsquigarrow D$ of the functor $\gM^*[\ivtd p]$ in Theorem~\ref{thm:Kisin}\eqref{thm:Kisin:equiv}. We set $D:=(\gM/u\gM)[\ivtd p]$ with $\vphi_D:=\vphi_\gM\bmod{u\gM}$, and define $\Fil^1 D_K $ as below:
%\begin{equation}
%\Fil^1 D_K := \frac{\langle\vphi_\gM(\gM)\rangle}{\langle\vphi_\gM(\gM)\rangle\cap\PP(u)\gM} \left[\ivtd p\right]\subset \frac{\gM}{\PP(u)\gM} \left[\ivtd p\right]\cong D_K  .
%\end{equation}
%(See \S1.2.7 and Lemmas~1.2.6, 1.3.13 of  \cite{kisin:fcrys} for more details.\footnote{Indeed, the construction in \cite{kisin:fcrys} was carried out  more generally; namely, for any $\gM\cong\gM^*(T)$ where $T\in\Replatsth$.}) By  \cite[Proposition~2.2.2]{kisin:fcrys}  this is the desired weakly admissible filtered $\vphi$-module
%\end{rmk}

Let us briefly discuss the torsion theory.

\begin{defn}\label{def:KisinTor}
For a non-negative integer $h$, a  \emph{torsion $\vphi$-module of height $\leqs h$} is a tuple $(\gM,\vphi_\gM)$, where $\gM$ is a finitely generated $p^\infty$-torsion $\Sig$-module with no non-zero $u$-torsion, and $\vphi_\gM:\gM\ra\gM$ is  a $\vphi$-semilinear endomorphism such that  $\PP(u)^h$ annihilates $\coker(1\otimes\vphi_\gM:\vphi^*\gM \ra \gM)$. Let $(\Mod/\Sig)^{\leqs h}$ denote the category of torsion $\vphi$-modules of height $\leqs h$ with the obvious notion of morphisms.
\end{defn}
Note that a nonzero $p^\infty$-torsion $\Sig$-module is of projective dimension $\leqs1$ if and only if it has no non-zero $u$-torsion.

%\begin{lem}\label{lem:KisinRaynaud}
%Any torsion $\vphi$-module $\gM$ of height $\leqs h$ fits into a short exact sequence
%\begin{equation}\label{eqn:KisinRaynaud}
%0\ra \wt\gM' \xra f \wt\gM \ra \gM \ra 0
%\end{equation}
%where $\wt\gM$ and $\wt\gM$ are $\vphi$-modules of height $\leqs h$ and all the arrows are $\vphi$-compatible.
%\end{lem} 
%\begin{proof}
%This is proved in \cite[Lemma~2.3.4]{kisin:fcrys} for $h=1$, but the argument generalizes to any $h$ with the following change: since $L:=\gM/\langle\vphi(\gM)\rangle$ is an $\Sig/(\PP(u)^h)$-module (instead of an $\fo_K$-module), we take $\wt L$ to be a finite free  $\Sig/(\PP(u)^h)$ which surjects onto $L$.
%\end{proof}  
\begin{lem}[{\cite[Lemma~2.3.4]{kisin:fcrys}}]\label{lem:KisinRaynaud}
Any   $\gM\in(\Mod/\Sig)^{\leqs1}$ fits into a short exact sequence
\begin{equation}\label{eqn:KisinRaynaud}
0\ra \wt\gM' \xra f \wt\gM \ra \gM \ra 0
\end{equation}
where $\wt\gM, \wt\gM'\in \phimodBT\Sig$ and all the arrows are $\vphi$-compatible.
\end{lem} 

%\begin{subequations}
To a torsion $\vphi$-module $\gM$ of height $\leqs h$, we can associate the following torsion $\Zp$-module with continuous $\GKinfty$-action:
\begin{equation}\label{eqn:TSigTors}
T^*_\Sig(\gM):=\Hom_{\Sig,\vphi}(\gM,\wh\Sig^{\ur}\otimes_{\Zp}\Qp/\Zp).
\end{equation}
One can easily check that  the exact sequence \eqref{eqn:KisinRaynaud} induces the following exact sequence:
\begin{equation}\label{eqn:TSigRedn}
0 \ra T^*_\Sig(\wt\gM) \xra{T^*_\Sig(f)} T^*_\Sig(\wt\gM')  \lra T^*_\Sig(\gM) \ra 0,
\end{equation}
where the surjective map is constructed as follows. Identifying $\wt\gM'[\ivtd p]$ with $\wt\gM[\ivtd p]$ by $f[\ivtd p]$, a $\Sig$-linear map $l:\wt\gM'\ra\wh\Sig^{\ur}$ can be viewed as $l:\wt\gM \ra \wh\Sig^{\ur}\otimes_{\Zp}\Qp$ and clearly $l \bmod {\wh\Sig^{\ur}}$ factors through $\gM\cong \wt\gM/\wt\gM'$. The arrow is defined by sending $l$ to $l \bmod {\wh\Sig^{\ur}}$. 
%\end{subequations}

\section{Review: Breuil modules}\label{sec:DefsBreuil}
We recall basic definitions and constructions to work with various $S$-modules introduced by Breuil. All results in this section are ``standard'' to experts, so we often omit the proof by giving references.

\subsection{The ring $S$}\label{subsec:setting}
Let $S$ be the $p$-adic completion of the divided power envelop of $W(k)[u]$ with respect to the ideal generated by $\PP(u)$. It can be shown that $S$ can be viewed as a subring of $ K _0[[u]]$ whose elements are precisely those of the form $\sum_{i\geq0} a_i \frac{u^i}{q(i)!}$, where $q(i):=\lfloor \frac{i}{e} \rfloor$ with $e:=\deg\PP(u)$, and $a_i\in W(k)$ converge to $0$ as $i\to\infty$.\footnote{The proof is by rearranging the terms of $\sum_{i\geqs0} b_i\frac{\PP(u)^i}{i!}$ and show that the coefficient $a_i$ of $\frac{u^i}{q(i)!}$ converges, for which boundedness of $\set{\frac{p^n}{n!}}$ suffices. (From $\lim_{i\to\infty}b_i=0$, we deduce $\lim_{i\to\infty}a_i=0$.) In particular, this argument works when $p=2$.} (In particular, $S$ is naturally a $\Sig$-algebra.) We define a differential operator $N:=-u\frac{d}{du}$ on $S$, and a ring endomorphism $\vphi:\self S$ which extends the Witt vector Frobenius map by $\vphi(u)=u^p$. We let $\Fil^1 S \subset S$ denote the ideal topologically generated by $\PP(u)^i/i!$ for $i\geqs 1$. Since $p|\vphi(f)$ for any $f\in\Fil^1S$, we can define a $\vphi$-semilinear map $\vphi_1:=\frac{\vphi}{p}:\Fil^1 S \ra S$. 

%Let us briefly record the setting of Proposition~\ref{prop:main}. Let $\gM\in\phimodBT\Sig$, and let $V$ be the unique crystalline $\GK$-representation extending $T^*_\Sig(\gM)[\ivtd p]$, and $T:=T^*_\Sig(\gM)$ viewed as a $\GKinfty$-stable $\Zp$-lattice in $V$. Put $D:=D^*_{\cris}(V)$, which necessarily satisfies $\gr^w D_K=0$ for $w\ne0,1$. %We will prove Proposition~\ref{prop:main} by constructing  from $\gM$ another \emph{$\GK$-stable} $\Zp$-lattice $T'\subset V$ via a strongly divisible $S$-lattice in $S\otimes_{W(k)}D$, and show that $T=T'$.

\subsection{Strongly divisible $S$-modules}
\label{subsec:D}
Let $D$ be an admissible filtered $\vphi$-module over $K$ such that $\gr^w D_K=0$ for $w\ne0,1$. Let us consider $\cD:=S\otimes_{W(k)}D$ equipped with the following structure:
\begin{enumerate} 
\item a $\vphi$-linear map $\vphi_{\cD}:=\vphi_S\otimes\vphi_D:\self{S\otimes_{W(k)}D} $,
\item a derivation $N_\cD:=N\otimes\id_D:\self{S\otimes_{W(k)}D} $ over $N:\self S$,
\item an $S[\ivtd p]$-submodule $\Fil^1\cD:=\set{m\in S\otimes_{W(k)}D|\ m\bmod{\PP(u)} \in\Fil^1D_K}$.
\end{enumerate}
See \cite{Breuil:GriffithsTransv} or \cite[\S2.2]{Breuil:IntegralPAdicHodgeThy} for more details (allowing higher Hodge-Tate weights).
We now define strongly divisible $S$-lattices in $\cD$, following \cite[Definition~2.2.1]{Breuil:IntegralPAdicHodgeThy} except that we do not impose any condition\footnote{I.e., we do not assume that $D$ has no non-zero weakly admissible quotient pure of slope $1$  as in \cite[\S2.1]{Breuil:IntegralPAdicHodgeThy}.} on $D$ when $p=2$:
\begin{defnsub}\label{def:BreuilMod}
A \emph{strongly divisible $S$-lattice in $\cD$} is an $S$-submodule $\M$ which satisfies the following properties:
\begin{enumerate}
\item $\M$ is a finite free $S$-module with $\M[\ivtd p]=\cD$
\item $\M$ is stable under $\vphi_\cD$ and $N$
\item we have $\vphi_{\cD}(\Fil^1\M)\subset p\M$ where $\Fil^1\M:=\M \cap \Fil^1\cD$.
\end{enumerate} 
\end{defnsub}
We put $\vphi_1:=\ivtd p(\vphi_{\cD}|_{\Fil^1\M}):\Fil^1\M \ra\M$, One can recover $\vphi_{\cD}|_{\M}$ from $\vphi_1$ as follows: 
\begin{equation*}
\vphi_{\cD}(m) = c_1\iv \vphi_1(\PP(u)m)
\end{equation*} 
for any $m\in\M$, where $c_1:=\vphi(\PP(u))/p$ which can be easily seen to be a unit. Note also that  $N_{\M}:=N_{\cD}|_{\M}$ and $\vphi_1$ satisfy the following  relation called ``Griffiths transversality'':
\begin{equation}\label{eqn;GriffithsTransv}
N_{\M}\circ\vphi_1 = c_1\iv\vphi_1\circ (\PP(u)\tim N_{\M}) :\Fil^1\M \ra \M
\end{equation}
We can axiomatize these structures and define \emph{strongly divisible $S$-modules of weight $\leqs1$} without mentioning the ambient module $\cD$. See \cite[\S2.2]{Breuil:IntegralPAdicHodgeThy} for more details.
%We say that $\M$ equipped with $\Fil^1\M\subset \M$, $\vphi_1:\Fil^1\M\ra\M$ and $N_{\M}$ is a \emph{strongly divisible $S$-module of weight $\leqs 1$} if it can be realized as a strongly divisible $S$-lattice in $S\otimes_{W(k)}D$ for some weakly admissible filtered $\vphi$-module $D$ with $\gr^w D_K = 0$ for any $w\ne0,1$.

\subsection{Construction of strongly divisible $S$-modules}\label{subsec:MM}
Let $\gM$ be an object of either $\phimodBT\Sig$ or $(\Mod/\Sig)^{\leqs1}$, and we will define some additional structure on the $S$-module  $\M:=S\otimes_{\vphi,\Sig}\gM$, which makes $\M$ a strongly divisible $S$-module if $\gM$ is $\Sig$-free. Using the $S$-linear map $1\otimes\vphi_\gM:\M\cong S\otimes_\Sig (\vphi^*\gM)\ra S\otimes_\Sig\gM$, we define an $S$-submodule  $\Fil^1\M\subset \M$ and $\vphi_1:\Fil^1\M \ra \M$ as follows:
\begin{subequations}\label{eqn:MM}
\begin{gather}
\label{eqn:Filh}\Fil^1\M:=\set{x\in\M|\ 1\otimes\vphi_\gM(x) \in \Fil^1S\otimes_\Sig \gM \subset S\otimes_\Sig \gM } \\
\label{eqn:vphir}\vphi_1:\Fil^1\M \xra{1\otimes \vphi_\gM} \Fil^1S\otimes_\Sig \gM \xra{\vphi_1\otimes1} S\otimes_{\vphi,\Sig}\gM = \M
\end{gather}
\end{subequations}
Let $\vphi_{\M}:=\vphi\otimes\vphi_\gM$, and note that $\vphi_1 = \ivtd p (\vphi_{\M}|_{\Fil^1\M})$. %One can also give an axiomatized definition of torsion $S$-modules $\M$ which should come from $\M\in(\Mod/\Sig)^{\leqs 1}$. 
\begin{exasub}\label{exa:MM}
Consider $\gM=\Sig$ equipped with $\vphi_\gM = \PP(u)\tim\vphi$. Then $\M = \Fil^1\M = S$ and $\vphi_1(1) = \vphi(\PP(u))/p$.
\end{exasub}
This functor $S\otimes_{\vphi,\Sig}(\cdot)$ is obviously faithful. We now claim that this functor is also exact. %The only non-trivial assertion to check is the following lemma, which can be directly verified from \eqref{eqn:Filh}:
\begin{lemsub}\label{lem:exactBreuil}
Assume that we have a $\vphi$-compatible short exact sequence $(\gM^\bullet):\ 0\ra\gM'\ra\gM\ra\gM''\ra 0$ where $\gM'$, $\gM$, and $\gM''$ are objects of either $(\Mod/\Sig)^{\leqs1}$ or $\phimodBT\Sig$. %one of the following holds:
% \begin{enumerate}
%\item $\gM', \gM, \gM''\in(\Mod/\Sig)^{\leqs1}$.
%\item $\gM', \gM, \gM''\in\phimodBT\Sig$.
%\item $\gM', \gM\in\phimodBT\Sig$, and  $\gM''\in(\Mod/\Sig)^{\leqs1}$.
%\end{enumerate}
%Set $\M:=S\otimes_{\vphi,\Sig}\gM$, and similarly define $\M'$ and $\M''$. Then the following sequence
%\begin{equation}\label{eqn:exactBreuil}
%0\ra\M'\ra\M\ra\M''\ra0
%\end{equation}
Then the sequence $S\otimes_{\vphi,\Sig}(\gM^\bullet)$ is short exact, respects $\vphi_1$'s, and induces a short exact sequence on $\Fil^1$'s.
\end{lemsub}
\begin{proof}
The only case when the exactness of $S\otimes_{\vphi,\Sig}(\gM^\bullet)$ is possibly non-trivial is when $\gM''\in(\Mod/\Sig)^{\leqs1}$. In this case, it is enough to show $\Tor^\Sig_1(\gM'',S)=0$. But since $\gM''$ is a successive extension of $\Sig/p\Sig$ as a $\Sig$-module, it is enough to show $\Tor^\Sig_1(\Sig/p\Sig,S)=0$, which is obvious. The rest of the claim can be directly checked from the definition \eqref{eqn:MM}.
\end{proof}

Suppose $\gM\in\phimodBT\Sig$. Set $D:= D^*_{\cris}(V_p(G'))$ for a $p$-divisible group $G'$ over $\fo_K$ which corresponds to $\gM$ up to isogeny by Theorem~\ref{thm:Kisin}(\ref{thm:Kisin:equiv}). We write $\cD:=S\otimes_{W(k)}D$ and $\M:=S\otimes_{\vphi,\Sig}\gM$.
\begin{propsub}\label{lem:N}
We use the notation as above. 
Then there exists a natural $S$-linear embedding $\eta:\M \hra \cD$ such that $\vphi_{\cD}\circ \eta = \eta\circ\vphi_{\M}$, $N_{\cD}(\eta(\M))\subset \eta(\M)$, and $\eta(\Fil^1\M) = \eta(\M)\cap\Fil^1\cD$.
\end{propsub}
\begin{proof}%[Sketch of the proof]
The assertion except $N_{\cD}(\M)\subset \M$ follows from \cite[Corollary~3.2.3]{Liu:StronglyDivLattice} using the rigid analytic construction \cite[Corollary~1.3.15]{kisin:fcrys}  relating $D$ and $\gM$.  
Now one can apply \cite[Proposition~5.1.3]{Breuil:GrPDivGrFiniModFil}\footnote{The proof of \cite[Proposition~5.1.3]{Breuil:GrPDivGrFiniModFil} works even when $p=2$.} to conclude.
\end{proof}
\begin{corsub}\label{cor:N}
For any $\gM\in(\Mod/\Sig)^{\leqs1}$, there exists a unique differential operator $N_{\M}:\M\ra\M$ over $N:S\ra S$ such that the ``Griffiths transversality'' \eqref{eqn;GriffithsTransv} holds and $N_{\M}\equiv 0 \bmod{I_0\M}$ where $I_0\subset S$ is the ideal topologically generated  by $\frac{u^{ei}}{i!}$ for $i\geqs0$.
\end{corsub}
\begin{proof}
By Lemma~\ref{lem:KisinRaynaud} we can find $\wt\gM,\wt\gM'\in\phimodBT\Sig$ with $\gM\cong \wt\gM/\wt\gM'$. By Proposition~\ref{lem:N}, both $\wt\M:=S\otimes_{\vphi,\Sig}\wt\gM$ and $\wt\M':=S\otimes_{\vphi,\Sig}\wt\gM'$ are strongly divisible $S$-lattices in $S\otimes_{W(k)}D$ for some filtered $\vphi$-module $D$. Therefore $N_{\M}:= N_{\wt\M}\bmod{\wt\M'}$ is well-defined and satisfies all the desired properties. The uniqueness follows from the proof of \cite[Lemma~3.2.1]{Breuil:IntegralPAdicHodgeThy}.
\end{proof}

\subsection{Crystalline period ring}
\label{par:Asthat}
Let us recall the construction of $\Acris$. For more complete discussions, see \cite[\S1, \S2]{fontaine:Asterisque223ExpII}.

There exists a ``canonical lift'' $\theta:W(\rep)\thra \fo_{\wh{\Kbar }}$ of the first projection $\rep \thra \fo_{\ol K }/(p)$, which is $\GK$-equivariant for the natural actions on both sides. % and is a topological quotient map (for the ``product topology'' on the source and the natural $p$-adic topology on the target).  %Let $\BdR^+$ be the completion of $W(\rep)[\ivtd p]$ with respect to the kernel of $\theta[\ivtd p]$, and let $\BdR:=\BdR^+[\ivtd t]$, where $t$ is the $p$-adic analogue of $2\pi i$. See \cite{fontaine:Asterisque223ExpII} for more details.
We define $\Acris$ as the $p$-adic completion of the divided power envelop of $W(\rep)$ with respect to $\ker(\theta)$. The Witt vector Frobenius map and the $\GK$-action on $W(\rep)$ extend to $\Acris$. We let $\Fil^1\Acris$ be the ideal topologically generated by $\ivtd{i!}(\ker\theta)^i$ for $i\geqs 1$ (under the $p$-adic topology). %The Witt vector Frobenius on $W(\rep)$ extends to $\Acris$, which we denote by $\vphi$. There exists a natural $W(\kbar)$-embedding  $\Acris\subset \BdR^+$ which respects $\vphi$, filtrations and $\GK$-actions.
We have $\vphi(\Fil^1\Acris)\Acris \subset p\Acris$, so we can define $\vphi_1:=\frac{\vphi}{p}:\Fil^1\Acris \ra \Acris$. 

The embedding $\Sig\hra W(\rep)$ constructed in \S\ref{subsec:GKinftyRep} extends to $S\hra \Acris$ by the universal property of divided powers envelop. This map is  $\vphi$-compatible and $\GKinfty$-invariant (but  \emph{not} $\GK$-invariant). 

Let us consider the sequence $\nf\epsilon:=\set{\epsilon^{(n)}}_{n\geqs0}$ of cocycle $\epsilon^{(n)}:\GK\ra\rep\starr$ defined as follows:
\begin{equation}
\epsilon^{(n)}(g):=  g\tim \pi^{(n)}/\pi^{(n)}
\end{equation}
for any $g\in\GK$. Clearly $\nf\epsilon(g)$ defines an element in $\rep$, and the $n$th component $\epsilon^{(n)}(g)$ is a (not necessarily primitive) $p^n$th root of unity. One can easily check that the following formula defines an element in $\Fil^1\Acris$:
\begin{equation}\label{eqn:tg}
t_g:=\log [\nf\epsilon(g)]=\sum_{n=1}^\infty\frac{(-1)^{n-1}([\nf\epsilon(g)]-1)^n}{n}
\end{equation}
where  $g\in\GK$ and $[\cdot]$ denotes the Teichm\"uller lift. %Conversely, we have $[\nf\epsilon(g)] = \sum_{n=0}^\infty \frac{t^n}{n!}$.

\begin{rmksub}\label{rmk:Acris}
When $p=2$, one can show that $t_g/2 \in\Fil^1\Acris$ (by modifying the argument in \cite[\S5.2.9]{fontaine:Asterisque223ExpII}), so we have $\gamma_n(t_g/2)\in\Fil^1\Acris$ since $\Fil^1\Acris$ has a divided power structure. It follows that for any $p$ the sequence $\gamma_n(t_g)$ converges  to $0$ as for the $p$-adic topology in $\Acris$. In particular, we have an equality $[\nf\epsilon(g)] = \sum_{n=0}^\infty \gamma_n(t_g) (=\exp(t_g))$ in $\Acris$. 
\end{rmksub}

Choose $g\in\GK$ such that $g\pi^{(1)}\ne\pi^{(1)}$. Set $[\nf\epsilon]:=[\nf\epsilon(g)]$ and $t:=t_g$ as in \eqref{eqn:tg}. For any integer $n\geqs0$ we set $ t^{\{n\}}:= t^n/(p^{q(n)}q(n)!)$ where $q(n):=\lfloor n/(p-1)\rfloor$. Note that $t^{\{n\}}\in\Acris$, and if $p=2$ then $ t^{\{n\}} = \gamma_n(t/2)$.% because  $t^{\{n\}} = t^{n-(p-1)q(n)}(t^{p-1}/p)^{q(n)}/q(n)!$. 

Consider the following subring $\fa_{\cris}$ of $\Frac W(\rep)[[T]]$, where $T$ is a formal variable:
\begin{equation*}
\fa_{\cris}:=\left\{\sum_{n=0}^\infty a_n T^{\{n\}} \,|  \ a_n\in W(\rep) \text{ such that } p\text{-adically } a_n\to 0. \right\}.
\end{equation*}
Now, one can give the following explicit description of $\Acris$.  
See \cite[ \S5.2.7--5.2.9]{fontaine:Asterisque223ExpII} for the proof.
\begin{thmsub}\label{thm:Acris}
The $W(\rep)$-algebra map $\fa_{\cris}\ra\Acris$ defined by $T^{\set n}\mapsto t^{\set n}$  is surjective and the kernel is generated by $[\nf\epsilon]-\sum_{n=0}^\infty \gamma_n(T)$.
\end{thmsub}
%
%We make the following interesting observation when $p=2$, which will be used later. Since $2$ divides $t_g^n/n!$ for any $n>0$, we have that $[\nf\epsilon(g)] \equiv 1 \bmod{2\Acris}$ for any $g\in\GK$. In particular, the subring $\Sig/(2)$ of $\Acris/(2)$ is fixed by the natural $\GK$-action on $\Acris/(2)$, so the subring $\wh\Sig^{\ur}/(2)$ of $\Acris/(2)$ is stable under the natural $\GK$-action on $\Acris/(2)$. 
%%By a formal identity, one obtains $[\nf\epsilon(g)]^p =  \sum_{n=0}^\infty (pt_g)^n/n! \equiv 1\bmod p\Acris$ for any $g\in\GK$, so  it follows that $\vphi([\nf\pi]) \equiv [\nf\epsilon(g)]^p[\nf\pi]^p\equiv \vphi([g\nf\pi]) \bmod p\Acris$. Therefore the image of $\vphi:\wh\Sig^{\ur}/(p)\ra\Acris/(p)$ is independent of the choice of $\nf\pi$ except $\pi=\pi^{(0)}$.
%
\subsection{Construction of Galois-stable $\Zp$-lattices}\label{par:Tst}
Let $\gM$ be an object of either $\phimodBT\Sig$ or $(\Mod/\Sig)^{\leqs1}$, and put $\M:=S\otimes_{\vphi,\Sig}\gM$ with additional structure defined in \S\ref{subsec:MM}. Consider the following $\Zp$-module to which we will endow a natural continuous $\GK$-action:
\begin{subequations}\label{eqn:Tst}
\begin{gather}
T^*_{\st}(\M):=\Hom_{S,\vphi_1, \Fil^1}(\M,\Acris) \text{, if }\gM\in\phimodBT\Sig, \\
T^*_{\st}(\M):=\Hom_{S,\vphi_1, \Fil^1}(\M,\Acris\otimes_{\Zp}\Qp/\Zp) \text{, if }\gM\in(\Mod/\Sig)^{\leqs1}.
\end{gather}
\end{subequations}
We let $\GKinfty$ act through the second argument $\Acris$. The $\GK$-action is more subtle to describe because the embedding $S\ra \Acris$ is only $\GKinfty$-stable, not $\GK$-stable. Let us first define a $\GK$-action on $\Acris\otimes_S\M$ using the differential operator $N_{\M}$ as follows.
\begin{subequations}\label{eqn:TstGK}
\begin{equation} \label{eqn:TstGK:a}
g\tim (a\otimes x) := g(a)\sum_{i=0}^\infty \gamma_i(t_g)\otimes N_{\M}^i(x), 
\end{equation}
for $g\in\GK$, $a\in\Acris$, and $ x\in\M$. Here, $t_g\in\Fil^1\Acris$ is as in \eqref{eqn:tg}, and $\gamma_i(t_g)$ is the standard $i$th divided power; i.e., $\gamma_i(t_g):=t_g^i/i!$ if $i>0$ and $\gamma_0(t_g):=1$ (even when $t_g=0$).  
%Here $\gamma_i$denotes the standard $i$th divided power on $\Fil^1\Acris$; i.e., $\gamma_i(a):= \frac{a^i}{i!}$ when $i>0$, and $\gamma_0(a):=1$ for any $a\in\Fil^1\Acris$ (including when $a=0$). 
By \cite[Lemma~5.1.1]{Liu:StronglyDivLattice}, the  sum (\ref{eqn:TstGK:a}) conerges and gives a $\GK$-action which respects $\vphi$ and the natural fintration on $\Acris\otimes_S\M$. Note also that we recover the natural $\GKinfty$-action on $\Acris\otimes_S\M$ since $t_g=0$ for $g\in\GKinfty$. %Since $\nf\epsilon$ is a $1$-cocycle, the formula \eqref{eqn:TstGK:a} defines a $\GK$-action, provided that it converges for any $g\in\GK$. 

For any $f\in T^*_{\st}(\M)$, we $\Acris$-linearly extend $f$ to $\Acris\otimes_S\M \ra \Acris$. Now it is clear that the following formula defines a continuous action of $g\in\GK$ on $T^*_{\st}(\M)$:
\begin{equation} \label{eqn:TstGK:b}
g\tim f: x\mapsto g\tim (f(g\iv(1\otimes x))), \text{ for }x\in\M.
\end{equation}
\end{subequations}
%\begin{rmksub}\label{rmk:GK}
%The following observation will be used in the proof of Proposition~\ref{prop:key}. Let $\gM\in\phimodBT\Sig$, and set $\M:=S\otimes_{\vphi,\Sig}\gM$. If $p=2$ then the $\GK$-action on $\Acris\otimes_S\M$ defined in \eqref{eqn:TstGK:a} simply becomes $g(a\otimes x) = g(a)\otimes x$ since  $\gamma_i(t_g)$ for $i>0$ are divisible by $2$ in $\Acris$ (\emph{cf.} Remark~\ref{rmk:Acris}). So the $\GK$-action on $T^*_{\st}(\M)$ defined in \eqref{eqn:TstGK:b} becomes $g\tim f:x\mapsto g\tim(f(x))$.
%\end{rmksub}

The following lemma is straightforward from \cite[Lemma~5.2.1]{Liu:StronglyDivLattice}\footnote{The assumption that $p>2$ under which \cite{Liu:StronglyDivLattice} was written is not relevant to this result.} and \cite[Proposition~2.2.5]{Breuil:IntegralPAdicHodgeThy}, using that $T^*_{\st}(\M)$ is $p$-adically separated and complete. %with no non-zero $p$-torsion.
\begin{lemsub}\label{lem:TstEmb}
%The formula \eqref{eqn:TstGK:b} defines a $\GK$-action on $T^*_{\st}(\M)$. Furthermore, if $\M$ is a strongly divisible $S$-lattice in $S\otimes_{W(k)}D$ then there exists a natural $\GK$-isomorphism $T^*_{\st}(\M)[\ivtd p] \cong V^*_{\cris}(D)$.
Let $D$ be as in \S\ref{subsec:setting}. If $\M$ is a strongly divisible $S$-lattice in $S\otimes_{W(k)}D$ then there exists a natural $\GK$-equivariant injective map $T^*_{\st}(\M) \hra V^*_{\cris}(D)$, whose image is a full $\Zp$-lattice.
\end{lemsub}

The following lemma on the exactness of $T^*_{\st}$  will be needed in \S\ref{sec:proofs}.
\begin{lem}\label{lem:TstExact}
Let $D_i$ for $i=1,2,3$ be admissible filtered $\vphi$-modules over $K$ with $\gr^w(D_i)_K=0$ for $w\ne 0,1$, and let $\M_i$ be a strongly divisible $S$-lattices in $S\otimes_{W(k)}D_i$.
For  any short exact sequence $0\ra\M_1\ra\M_2\ra\M_3\ra0$, the following sequence is short exact:
\begin{equation}\label{eqn:TstExact}
0 \ra T^*_{\st}(\M_3) \ra T^*_{\st}(\M_2) \ra T^*_{\st}(\M_1) \ra 0. 
\end{equation} 
%If $\gM_1,\gM_2\in\phimodBT\Sig$ and $\gM_3\in(\Mod/\Sig)^{\leqs1}$ then the following sequence
%\begin{equation}\label{eqn:TstExact:Redn}
%0 \ra T^*_{\st}(\M_2) \ra T^*_{\st}(\M_1) \ra T^*_{\st}(\M_3) \ra 0, 
%\end{equation} 
%defined in the similar manner to \ref{eqn:TSigRedn}, is short exact.
\end{lem}
\begin{proof}
The left exactness of (\ref{eqn:TstExact}) is clear from the definition of $T^*_{\st}$, and the right exactness can be check after reducing it modulo $p$. Moreover, it suffices to show the left exactness of the mod~$p$ reduction of  (\ref{eqn:TstExact}), as the $\Fp$-dimension of each term adds up correctly.

Now, the natural map $T^*_{\st}(\M_i)/(p)\ra T^*_{\st}(\M_i/(p))$ is injective, and the natural map $T^*_{\st}(\M_3/(p)) \ra T^*_{\st}(\M_2/(p))$ induced by $\M_2/(p)\thra \M_3/(p)$ is injective as $T^*_{\st}$ is  left exact in the torsion case as well. This shows the left exactness  the mod~$p$  reduction of (\ref{eqn:TstExact}), hence the lemma follows.
\end{proof}

\subsection{Relation with Galois representations coming from Kisin modules}\label{subsec:imath}
%Let us first construct an embedding $\wh\Sig^{\ur}\hra\Acris$ of $\Sig$-algebras which commutes with $\vphi$ and $\GKinfty$-action. Wintenberger constructed a continuous $\GKinfty$-equivariant embedding $\Lambda_{\Kbar/K_\infty/K}:\ol X \hra \Frac\rep$ over $X_K(K_\infty)=k\llpar u \rrpar$ with dense image \cite[Corollaire~4.3.4]{wintenberger:NormFiels}
%, which unique up to $\GKinfty$-action (\emph{cf.} \S\ref{par:Asthat}). By the universal property of strict henselization, there exists a unique embedding $\imath: \wh\fo_\Eps^{\ur}\hra W(\Frac\rep)$ of $\fo_\Eps$-algebras which lifts the chosen embedding $k\llpar u \rrpar^{\sep} \hra \Frac\rep$. It follows (again from the universal property of strict henselization) that  $\imath$ commutes with $\vphi$ and $\GKinfty$-action. 
%
%Since $\imath(\Sig)\subset W(\rep)$ and $W(\rep)$ is integrally closed and topologically closed in $W(\Frac\rep) $, it follows that $\imath$ restricts to $\imath:\wh\Sig^{\ur}\hra W(\rep)$, and hence we obtain $\wh\Sig^{\ur}\hra\Acris$. Note that the choice we make modifies the map only by $\GKinfty$-action.
Note that the natural map $W(\rep)\ra\Acris$ is injective \cite[\S2.3.3]{fontaine:Asterisque223ExpII}, so we may view $\wh\Sig^{\ur}$ as a subring of $\Acris$ (\emph{cf.} \S\ref{subsec:GKinftyRep}). Now consider the following natural $\GKinfty$-equivariant morphism: 
\begin{equation}\label{eqn:TSigTqstLatt}
\imath: T^*_\Sig(\gM) = \Hom_{\Sig,\vphi}(\gM,\wh\Sig^{\ur}) \ra \Hom_{S,\vphi_1,\Fil^1}(\M,\Acris) = T^*_{\st}(\M)
\end{equation}
where the arrow in the middle is defined as follows: the arrow in \eqref{eqn:TSigTqstLatt} takes $f:\gM\ra\wh\Sig^{\ur}$ to $\wt f:\M = S\otimes_{\vphi,\Sig}\gM \ra \Acris$ which is obtained by $S$-linearly extending $\gM \xra f \wh\Sig^{\ur}\xra\vphi \Acris$, where we view $S$ as a $\Sig$-algebra via $\vphi:\Sig\ra S$. One can check that if $f$ respects $\vphi$, then $\wt f$ respects $\vphi_1$ and takes $\Fil^1\M$ to $\Fil^1\Acris$. This map \eqref{eqn:TSigTqstLatt} is clearly $\GKinfty$-equivariant, and  one can direct check that the map  \eqref{eqn:TSigTqstLatt} respects the natural embeddings of the source and the target into $V$.\footnote{See the proof of \cite[Proposition~2.1.5]{kisin:fcrys} for $ T^*_\Sig(\gM)\hra V^*_{\cris}(D)$, and the proofs of \cite[Proposition~2.2.5]{Breuil:IntegralPAdicHodgeThy} and  \cite[Lemma~5.2.1]{Liu:StronglyDivLattice} for $ T^*_{\st}(\M)\hra V^*_{\cris}(D)$.}

%Let $T'=T^*_{\qst}(\M)$ denote the $\GKinfty$-stable $\Zp$-lattice in $V$, which is also $\GK$-stable by Lemma~\ref{lem:QstVsSt}. 
Unfortunately, when $p=2$ the map $\imath$ \eqref{eqn:TSigTqstLatt} is \emph{not} in general an isomorphism.\footnote{If $p>2$ then map $\imath$ is an isomorphism. The proof can be found in the proof of \cite[Theorem~2.2.7]{kisin:fcrys}, and it can also be deduced from Proposition~\ref{prop:main2} below.} We prove, instead, that the image of $\imath$ is a $\GK$-stable submodule in $T^*_{\st}(\M)$, which suffices for our purpose (\emph{cf.} Proposition~\ref{prop:main}). This claim will proved in \S\ref{sec:proofs}.%But a slightly weaker statement is true which is enough for proving Proposition~\ref{prop:main}. 
%\begin{prop}\label{prop:mainRestated}
%The image of $\imath$ in $T^*_{\st}(\M)$ is $\GK$-stable.
%\end{prop}
%We will prove the following proposition in the next section.
%
%\section{Kisin modules and strongly divisible $S$-lattices}\label{sec:IntHT}

\section{Statement of the main result}\label{sec:statements}
Let $G$ be a $p$-divisible group over $\fo_K$, and $T:=T_p(G)$ its Tate module. %, and $V:=T[\ivtd p]$. Since $V$ is a crystalline $\GK$-representation with  $\gr^wD_{\dR}^*(V) = 0$ for $w\ne 0,1$ \cite[Theorem~6.2]{Fontaine:BT}, we can use  Kisin theory (Theorem~\ref{thm:Kisin}\eqref{thm:Kisin:main}) to  obtain the following functor $\gM^*\circ T_p$:
Recall that we have the following functor defined in Theorem~\ref{thm:Kisin}\eqref{thm:Kisin:equiv}:
\begin{equation*}\label{eqn:gM}
%\gM^*\circ T_p :
\left\{p\text{-divisible groups over }\fo_K\right\} \ra \phimodBT\Sig;\ G \rightsquigarrow \gM_G:=\gM^*(T_p(G)).
\end{equation*}
By the theorems of Tate \cite[(4.2)]{Tate:pDivGps} and Kisin (Theorem~\ref{thm:Kisin}\eqref{thm:Kisin:main}), the functor $\gM^*\circ T_p$ is fully faithful. %Furthermore, we will show the following theorem:

\begin{thm}\label{thm:classif}
Assume that $\fo_K$ is a $2$-adic discrete valuation ring with perfect residue field.
Then the functor $\gM^*\circ T_p$  is an anti-equivalence of categories. Furthermore, for any $p$-divisible group $G$  over $\fo_K$ and $\gM_G:=\gM^*(T_p(G))$, there is a natural $\GKinfty$-isomorphism $T_p(G)\cong T^*_\Sig(\gM_G)$.
\end{thm}
Note that the second assertion is automatic from the property of $\gM^*$ stated in Theorem~\ref{thm:Kisin}\eqref{thm:Kisin:main}. 

Now, let us state the relation of this classification with crystalline Dieudonn\'e theory. Let $G$ be a $p$-divisible group over $\fo_K$ with residue characteristic $p=2$, $\gM_G:=\gM^*(T_p(G))$, and $D:=D^*_{\cris}(V_p(G))$. Let $\DD^*(G)$ denote the contravariant Dieudonn\'e crystal. (See, for example, \cite{Mazur-Messing} for the definition.) 
Then, $\DD^*(G)(S)$ and $S\otimes_{\vphi,\Sig}\gM_G$ are strongly divisible $S$-lattices in $S\otimes_{W(k)}D$. (See \S\ref{subsec:MM} and  \cite[\S{A}]{kisin:fcrys} for more details.) We will prove the following proposition in \S\ref{subsec:Dieudonne} and \S\ref{subsec:exactness}.
\begin{prop}\label{prop:Dieudonne}
We use the same notation as above.
\begin{enumerate} 
\item
\label{prop:Dieudonne:Comp}
We have $\DD^*(G)(S) = S\otimes_{\vphi,\Sig}\gM_G$ as strongly divisible $S$-lattices in $S\otimes_{W(k)}D$.
\item
\label{prop:Dieudonne:Exact}
The functor $\gM^*\circ T_p$ takes a short exact sequence of $p$-divisible groups to a short exact sequence in $\phimodBT\Sig$.
\end{enumerate}
\end{prop}
The proof of \eqref{prop:Dieudonne:Exact} uses \eqref{prop:Dieudonne:Comp} and the exactness  of $\DD^*$; as explained in Remark~\ref{rmk:PathologyLiu} we cannot hope to prove the exactness of $\gM^*\circ T_p$ by a purely linear algebraic method. We also prove a stronger exactness assertion which is useful in practice; namely, that any quasi-inverse of $\gM^*\circ T_p$ should also be an exact functor.
%
%Both Theorem~\ref{thm:classif} and Proposition~\ref{prop:Dieudonne} hold when $p>2$, which were proved by Kisin \cite[\S2.2, \S{A}]{kisin:fcrys}. When $p=2$, Kisin  \cite[\S1]{Kisin:2adicBT} proved them for  the Cartier duals of connected $p$-divisible groups, and Lau \cite{Lau:2010fk} constructed an exact  (anti-)equivalence\footnote{Lau \cite{Lau:2010fk} worked with covariant Dieudonn\'e crystals, so his convention differs from ours by Cartier duality.} of categories between the category of $p$-divisible groups over $\fo_K$ and $\phimodBT\Sig$ which satisfies Proposition~\ref{prop:Dieudonne}, but he did not showed the compatibility of his (anti-)equivalence with ``Kisin theory'' for the Tate modules (hence he did not obtain the second assertion of Theorem~\ref{thm:classif}). 

Raynaud's theorem \cite[Th\'eor\`eme~3.1.1]{Berthelot-Breen-Messing:DieudonneII} implies that for any $p$-power order finite flat group scheme $H$ over $\fo_K$, there exists an isogeny $f:G \ra G'$ of $p$-divisible groups over $\fo_K$ such that $H\cong \ker (f)$. We set $\gM_H:=\coker\left((\gM^*\circ T_p)(f)\right)$, and one can check that $\gM_H\in  (\Mod/\Sig)^{\leqs1}$ does not depend on the choice of $f$ up to isomorphism, using Proposition~\ref{prop:Dieudonne}(\ref{prop:Dieudonne:Exact}). One can also define the (contravariant) Dieudonn\'e crystal $\DD^*(H)$ in such a way that we have the following natural short exact sequence
\begin{equation*}
0 \ra \DD^*(G')(S) \xra{\DD^*(f)} \DD^*(G)(S) \ra \DD^*(H)(S) \ra 0
\end{equation*}
 (See D\'efinition~3.1.5 and Th\'eor\`eme~3.1.2 of \cite{Berthelot-Breen-Messing:DieudonneII} for more details.) From this exact sequence one can define $\Fil^1$ and $\vphi_1$ on  $\DD^*(H)(S)$,\footnote{I.e.,  $\DD^*(H)(S)\in(\Mod/S)$ in the sense of \cite[\S2.1.1]{Breuil:GrPDivGrFiniModFil}, but allowing $p=2$.} which do not depend on the choice of $f$.
\begin{cor}\label{cor:classifFF}
There exists an exact anti-equivalence of categories $H\rightsquigarrow \gM_H$ from the category of $p$-power order finite flat group schemes over $\fo_K$ to $(\Mod/\Sig)^{\leqs1}$ which satisfies the following properties:
\begin{enumerate}
\item There is a natural $\GKinfty$-isomorphism $H(\Kbar) \cong T^*_\Sig(\gM_H)$, where $T^*_\Sig$ is defined in \eqref{eqn:TSigTors}.
\item For any $p$-divisible group $G$ over $\fo_K$, we have a natural isomorphism $\gM^*(T_p(G)) \cong \varprojlim_n \gM_{G[p^n]}$, where the projective limit is taken with respect to the map induced by the natural inclusions $G[p^n]\hra G[p^{n+1}]$.
\item There exists a natural isomorphism $S\otimes_{\vphi,\Sig}\gM_H \cong \DD^*(H)(S)$ which respects $\Fil^1$ and $\vphi_1$.
\end{enumerate}
\end{cor}
\begin{proof}
Following the same proof of \cite[Theorem~2.3.5]{kisin:fcrys} one can deduce the corollary from Theorem~\ref{thm:classif}, Proposition~\ref{prop:Dieudonne}, and Lemma~\ref{lem:KisinRaynaud} (with the discussion following it), possibly except the exactness of the functor $H\rightsquigarrow \gM_H$ which is proved in Proposition~\ref{prop:exact}. Note that the argument crucially uses the \emph{exactness} of the functor $\gM\circ T_p$.
\end{proof}

This concludes the proof of Conjecture~\ref{conj:breuil}. See the introduction (\S\ref{sec:intro}) for previous results on this conjecture.

It was observed by Breuil \cite[Theorem~3.4.3]{Breuil:IntegralPAdicHodgeThy} that  the following corollary is a formal consequence\footnote{Since we have the classification of any $p$-power order finite flat group schemes, not just the ones killed by $p$, the d\'evissage step in the proof of \cite[Theorem~3.4.3]{Breuil:IntegralPAdicHodgeThy} is unnecessary.} of Corollary~\ref{cor:classifFF}. We refer to \emph{loc.cit} for the proof.
\begin{cor}\label{cor:FF}
The functor  $T\rightsquigarrow T|_{\GKinfty}$ from the category of $p^\infty$-torsion $\GK$-representations coming from some finite flat group scheme over $\fo_K$ to the category of $p^\infty$-torsion $\GKinfty$-representations is fully faithful. 
\end{cor}

Let us sketch the strategy of the proof of Theorem~\ref{thm:classif}, which can be thought of as extending \cite[\S2.2]{kisin:fcrys} to the case when $p=2$. The only assertion that needs to be proved is the essential surjectivity of $\gM^*\circ T_p$; i.e. given any $\gM\in\phimodBT\Sig$, we need to produce a $p$-divisible group $G$ over $\fo_K$ such that $\gM\cong \gM^*(T_p(G))$. For any given $\gM\in\phimodBT\Sig$,  Kisin's theorem (Theorem~\ref{thm:Kisin}\eqref{thm:Kisin:equiv}) provides a $p$-divisible group $G'$ over $\fo_K$ such that $\gM[\ivtd p]\cong \gM^*(T_p(G'))[\ivtd p]$. Set $V:= T_p(G')[\ivtd p]$ and note that $T^*_\Sig(\gM)$ embeds into $V$ as a $\GKinfty$-stable $\Zp$-lattice.
The main step is to show the following proposition:

\begin{prop}\label{prop:main}
With the notation as above, the $\GKinfty$-stable $\Zp$-lattice $T$ in $V$ is $\GK$-stable. 
\end{prop}
Granting this proposition, $T$ comes from a $p$-divisible group over $\fo_K$ by Raynaud's ``limit theorem'' \cite[Proposition~2.3.1]{raynaud:GpSch}, so we obtain Theorem~\ref{thm:classif}. 

The proof of Proposition~\ref{prop:main} (which will be given in \S\ref{sec:proofs}) is not straightforward since Breuil's theory of strongly divisible $S$-lattices does not work nicely when $p=2$, but fortunately it does not fail too badly. (See Propositions~\ref{prop:main2}, \ref{prop:key} for more details.) Interestingly, we do not use Proposition~\ref{prop:Dieudonne} in order to prove Proposition~\ref{prop:main} or Theorem~\ref{thm:classif}.
%
%Finally, the proof of Proposition~\ref{prop:Dieudonne} could be reduced (by some trick explained in \S\ref{subsec:Dieudonne}) to the cases when the Cartier dual $G^\vee$ is either connected or \'etale, where the theorem is already known by Kisin \cite[Proposition~1.1.10]{Kisin:2adicBT}. 

\begin{rmk}\label{rmk:IntModel}
Kisin's construction of integral canonical models for Shimura varieties of abelian type \cite{Kisin:IntModelAbType} has a few technical restrictions when $p=2$, and one of the main reasons is that at the time of \cite{Kisin:IntModelAbType} we did not have Theorem~\ref{thm:classif} and Proposition~\ref{prop:Dieudonne} when $p=2$. (See \cite[Theorem~1.4.2, Corollary~1.4.3]{Kisin:IntModelAbType} for more details.) The $2$-adic case causes another complication (of different nature) when the reductive group $G$ has a factor of type~B \cite[Lemma~2.3.1]{Kisin:IntModelAbType}\footnote{In concrete cases, this type~B restriction is not so serious; this restriction arises when finding a suitable embedding $G\hra GL_n$ over $\Zp$, which is often provided by definition for most groups $G$ that arise in practice.}, % For example, we often get a reductive group $G$ of type~B from a certain odd dimensional quadratic space $V$, which naturally comes equipped with a ``good'' embedding  $G\hra \GL(V)$.}, 
so one may expect that the main result of this paper can be used to construct $2$-adic integral models for Shimura varieties of abelian type under some mild restriction on the reductive group $G$. Unfortunately, there is one more deformation argument \cite[Proposition~1.5.8]{Kisin:IntModelAbType} which does not generalize to the $2$-adic case. This problem, however, does not seem insurmountable, and it seems that the main result of this paper has removed the main obstacle for generalizing \cite{Kisin:IntModelAbType} to the $2$-adic case. 
\end{rmk}

\section{Proof of the main results}\label{sec:proofs}
In this section, %we prove Proposition~\ref{prop:mainRestated} (hence, Theorem~\ref{thm:classif}) and Proposition~\ref{prop:Dieudonne}. 
we prove  Proposition~\ref{prop:main} (hence, Theorem~\ref{thm:classif}) and Proposition~\ref{prop:Dieudonne}. The result is known when $p>2$ by Kisin \cite[\S2.2]{kisin:fcrys} or when $\gM$ is $\vphi$-unipotent (\emph{cf.} \S\ref{par:DefEtPhiNilp}) by Kisin \cite[Proposition~1.2.7]{Kisin:2adicBT}, and our proof can be viewed as a modification of these arguments. A similar approach to ours (namely, classifying connected and \'etale $p$-divisible groups separately and handling the general case by connected-\'etale sequences) is also found in \cite{BlockKato:Dieudonne} and \cite{dejong:crysdieubyformalrigid} in their study of Dieudonn\'e crystals, and  in \cite{Zink:DieudonnePDivGpCFT} where display theory is extended to general $p$-divisible groups.
%Let us state the key lemma, which crucially uses all of the basic assumptions in \S\ref{par:BasicAssumpBreuilMod}. We postpone the proof to  \S\ref{par:TongLemmaPf}.
\subsection{}\label{par:DefEtPhiNilp}
We first need to introduce the linear algebra counterparts of \'etale, connected, multiplicative-type, and unipotent finite flat group schemes and $p$-divisible groups over $\fo_K$, respectively. %, for $\vphi$-module of height $\leqs 1$. 
Let $\gM$ be an object either in $(\Mod/\Sig)^{\leqs1}$ or in $\phimodBT\Sig$. Since the linearlzation $1\otimes\vphi_\gM:\vphi^*\gM\ra\gM$ is injective with cokernel killed by $\PP(u)$, there exists a unique map $\psi_\gM:\gM\ra\vphi^*\gM$ such that $\psi_\gM\circ(1\otimes\vphi_\gM)$ and $(1\otimes\vphi_\gM)\circ\psi_\gM$ are given by multiplication by $\PP(u)$. 
 
 We first define the analogue of Cartier duality $\gM^\vee$ as follows. We set the underlying $\Sig$-module to be $\gM^\vee:=\Hom_\Sig(\gM,\Sig)$ if $\gM\in\phimodBT\Sig$, and $\gM^\vee:=\Hom_\Sig(\gM, \Sig\otimes_{\Zp} \Qp/\Zp)$ if $\gM\in(\Mod/\Sig)^{\leqs1}$, and we define $\vphi_{\gM^\vee}$ so that  $1\otimes\vphi_{\gM^\vee}$ is the dual of $\psi_\gM$. One can check that $T^*_\Sig(\gM^\vee)\cong (T^*_\Sig(\gM))^*\otimes(\chi_{\cyc}|_{\GKinfty})$, so for  a $p$-divisible group $G$ over $\fo_K$ we have $\gM^*(T_p(G))^\vee\cong  \gM^*(T_p(G^\vee))$ where $G^\vee$ is the Cartier dual of $G$. (This essentially follows from $\gM^*(\Zp(1)) \cong (\Sig,\PP(u)\vphi)$; indeed, by \cite[Corollary~1.3.15]{kisin:fcrys} the isogeny class of $(\Sig,\PP(u)\vphi)$ corresponds to the admissible filtered $\vphi$-module $(W(k)[\ivtd p], \PP(0)\vphi)$, and from our normalization we have  $\PP(0) = p$.) The analogous statement for finite flat group schemes is deduced from this using Proposition~\ref{prop:Dieudonne}(\ref{prop:Dieudonne:Exact}) which will be proved later.%\footnote{We need the exactness of $\gM^*\circ T_p$ to show that $\gM_H$ does not depend on the embedding of $H$ into a $p$-divisible group.} 
 %(which we will prove later in \S\ref{subsec:exactness}). If this is done, then the claim follows because for any short exact sequence $0\ra\wt\gM'\ra\wt\gM\ra\gM\ra0$ with $\wt\gM,\wt\gM'\in\phimodBT\Sig$ and $\gM\in(\Mod/Sig)^{\leqs1}$ we naturally obtain a short exact sequence $0\ra \wt\gM^\vee\ra(\wt\gM')^\vee\ra\gM^\vee\ra0$.}
 
We say $\gM$ is \emph{\'etale}  if $1\otimes\vphi_{\gM}$ is an isomorphism, and \emph{of multiplicative type} if $\gM^\vee$ is \'etale (i.e., if $\psi_\gM$ is an isomorphism). We say $\gM$ is \emph{$\vphi$-nilpotent}  if $\vphi^n_{\gM}(\gM)\subset\m_\Sig\tim \gM$ for some $n\gg1$, and \emph{$\vphi$-unipotent}\footnote{Note that $\gM$ being $\vphi$-unipotent does \emph{not} imply that there exists a unipotent matrix representation of $\vphi_{\gM}$. This notion should be thought of as an analogue of unipotent $p$-divisible groups/finite flat group schemes.} if $\gM^\vee$ is $\vphi$-nilpotent (i.e., if $\psi^n_{\gM}(\gM)\subset\m_\Sig\tim (\vphi^n)^*\gM$ for some $n\gg1$). 

The following lemma asserts the existence of the linear algebra counterpart of connected-\'etale and multiplicative-unipotent sequences:
\begin{lem}\label{lem:ConndEtSeq}
%\begin{subequations}
For any $\gM\in\phimodBT\Sig$, there exist functorial short exact sequences
\begin{eqnarray}
\label{eqn:ConndEtSeq}0 \ra \gM^{\et}\ra\gM\ra\gM^c\ra0\\
\label{eqn:ConndEtDual} 0 \ra \gM^{u}\ra\gM\ra\gM^m\ra0
\end{eqnarray}
where $\gM^{\et},\gM^{u}\in\phimodBT\Sig$ are respectively maximal among \'etale and $\vphi$-unipotent submodules of $\gM$,  and $\gM^c, \gM^{m}\in\phimodBT\Sig$ are respectively maximal among $\vphi$-nilpotent and multiplicative-type quotients of $\gM$. Furthermore, we have $(\gM^{\et})^\vee = (\gM^\vee)^m$ and $(\gM^c)^\vee = (\gM^\vee)^u$. The same statement holds for $\gM\in(\Mod/\Sig)^{\leqs1}$.
%\begin{equation}\label{eqn:ConndEtSeq}
%0 \ra \gM^{\et}\ra\gM\ra\gM^c\ra0,
%\end{equation}
%where $\gM^{\et}\in\phimodBT\Sig$ is maximal among $\vphi$-stable \'etale submodule of $\gM$, and $\gM^{c}\in\phimodBT\Sig$ is maximal among $\vphi$-nilpotent quotient of $\gM$. The same holds for $\gM\in(\Mod/\Sig)^{\leqs1}$.
\end{lem}
\begin{proof}
When $\gM\in\phimodBT\Sig$, the lemma directly follows from \cite[Proposition~1.2.11]{Kisin:ModuliGpSch} applied to $\gM/(p^n)$ (which is free over $\Sig\otimes_{\Zp}\Zp/(p^n)$). When $\gM\in(\Mod/\Sig)^{\leqs1}$, choose $\gM\cong\wt\gM/\wt\gM'$ with $\wt\gM,\wt\gM'\in\phimodBT\Sig$ and one can check that the sequences induced from (\ref{eqn:ConndEtSeq}, \ref{eqn:ConndEtDual}) for $\wt\gM$ work.
\end{proof}

%By ``Cartier duality'' introduced in \S\ref{par:DefEtPhiNilp}, we obtain the following functorial short exact sequence for any $\gM\in\phimodBT\Sig$.
%\begin{equation}\label{eqn:ConndEtDual}
%0 \ra \gM^{u}\ra\gM\ra\gM^m\ra0
%\end{equation}
%where $\gM^{u}\in\phimodBT\Sig$ is maximal among $\vphi$-unipotent submodule of $\gM$,  and $\gM^{m}\in\phimodBT\Sig$ is the maximal among multiplicative-type quotient of $\gM$. (The same statement holds for $\gM\in(\Mod/\Sig)^{\leqs1}$.)
%%\end{subequations}
%
\begin{prop}\label{prop:Connd}
For a $p$-divisible group $G$ over $\fo_K$, put $\gM_G:=\gM^*(T_p(G))$.
Then the functor $\gM^*\circ T_p$ takes the connected-\'etale and multiplicative-unipotent sequences for $G$ to the exact sequences \eqref{eqn:ConndEtSeq} and \eqref{eqn:ConndEtDual} for $\gM_G$, respectively.
%The same statement holds for a $p$-power order finite flat group scheme $H$ over $\fo_K$ if $H\rightsquigarrow \gM_H$ is well-defined.
\end{prop}
%\begin{corsub}
%A $p$-divisible group is \'etale (respectively, connected) if and only if $\gM_G$ is \'etale (respectively, $\vphi$-nilpotent). The same statement holds for a $p$-power order finite flat group scheme $H$ over $\fo_K$ if $H\rightsquigarrow \gM_H$ is well-defined.
%\end{corsub}
We will show later (\S\ref{subsec:exactness})  that the functor $\gM^*\circ T_p$ takes any exact sequence of $p$-divisible groups to an exact sequence in $\phimodBT\Sig$ (not just connected-\'etale sequences), and its proof uses the compatibility with crystalline Dieudonn\'e theory whose proof uses Proposition~\ref{prop:Connd}.

We first need the following lemma which shows an analogue to the fact that any \'etale group scheme over $\fo_{\wh K^{\ur}}$ is constant.
\begin{lemsub}\label{lem:classifMult}
Assume that the residue field $k$ of $\fo_K$ is algebraically closed. Then $\gM\in\phimodBT\Sig$ is \'etale if and only if $\gM$  is isomorphic to the product of  copies of $(\Sig,\vphi)$. Also, $\gM\in\phimodBT\Sig$ is of multiplicative type if and only if $\gM$ is isomorphic to the product of copies of $(\Sig, \PP(u)\vphi)$.
\end{lemsub}
\begin{proof}
Since the two assertions are Cartier dual to each other (\emph{cf.} \S\ref{par:DefEtPhiNilp}), it is enough to prove the \'etale case. The ``if'' direction is obvious, so we focus on the ``only if'' direction. Assume that $\gM\in\phimodBT\Sig$ is \'etale, and we fix a $\Sig$-basis $\set{\e_1,\cdots,\e_n}$ of $\gM$. Let $A$ denote the matrix representation of $\vphi_\gM$ with respect to the chosen basis. If we replace this basis by $\set{U\e_i}$ for some invertible matrix $U$, then $A$ gets replaced by $\vphi(U)AU\iv$.

We first show that one can choose an invertible matrix $U_0$ over $\Sig$ such that $\vphi(U_0)AU_0\iv$ is congruent to $1$ modulo $(u)$; indeed, finding an $n\times n$ invertible matrix $\bar U_0$ over $W(k)$ such that $\vphi(\bar U_0) (A\bmod u)\bar U_0\iv = 1$ boils down to solving equations all of whose roots lie in some unramified extension of $W(k)$. So we may assume that $A=1+uB$ for some matrix $B$. Now, by taking $U_1=A$, we obtain $A_1:=\vphi(A)AA\iv = \vphi(A) = 1+u^p\vphi(B)$. Recursively, we take $U_{i+1}:=A_{i}$ and $A_{i+1} = \vphi(U_{i+1})A_i U_{i+1}\iv$. Since $U_{i+1}:=A_{i} \equiv 1 \bmod{u^{p^i}}$, the infinite product $U:=\prod_{i\geqs1} U_i$ converges and $\vphi(U)AU\iv = 1$.
\end{proof}

\begin{lemsub}\label{lem:etqt}
Let $\gM\in\phimodBT\Sig$, and let $\gM^{\et}$ be the maximal \'etale submodule of $\gM$ which exists by Lemma~\ref{lem:ConndEtSeq}. Then, $T^*_\Sig(\gM^{\et})$ is the maximal unramified $\GKinfty$-quotient of $T^*_\Sig(\gM)$.
\end{lemsub}
\begin{proof}
By Lemma~\ref{lem:classifMult}, $T^*_\Sig(\gM^{\et})$ is an unramified quotient of  $T^*_\Sig(\gM)$. To show it is maximal, it is enough to show that $\gM^{\et}\ne 0$  if $T^*_\Sig(\gM)$ has a non-zero unramified quotient. For this, choose  $f\in T^*_\Sig(\gM) = \Hom_{\Sig,\vphi}(\gM,\wh\Sig^{\ur})$ whose image in the maximal unramified quotient is not zero. Since $(\wh\Sig^{\ur})^{I_K} = W(\kbar)[[u]]$, it follows that $f\iv (W(\kbar)[[u]])$ is a non-zero \'etale submodule of $\gM$.
\end{proof}

\begin{proof}[Proof of Proposition~\ref{prop:Connd}]
%Note that Proposition~\ref{prop:Connd}(\ref{prop:Connd:et}) follows from Proposition~\ref{prop:Connd}(\ref{prop:Connd:seq})
Since two assertions in the statement are Cartier dual to each other (\emph{cf.} \S\ref{par:DefEtPhiNilp}), it suffices to show that  $\gM^*\circ T_p$ respects connected-\'etale sequences.
Let $(G^\bullet):\, 0\ra G^c\ra G\ra G^{\et}\ra 0$ be a connected-\'etale sequence of $p$-divisible groups, and $(\gM^\bullet):\, 0\ra\gM^{\et}_G\ra \gM_G \ra \gM_G^c \ra 0$ denote the ``connected-\'etale sequence'' for $\gM_G:=\gM^*(G)$ as in Lemma~\ref{lem:ConndEtSeq}. To prove Proposition~\ref{prop:Connd}, it suffices to show that $T_p(G^\bullet) \cong T^*_\Sig(\gM^\bullet)$ via the natural isomorphism. 

By Lemma~\ref{lem:etqt}, $T^*_\Sig(\gM_G^{\et})$ is the maximal unramified $\GKinfty$-quotient of $T_p(G)$. Viewing $T^*_\Sig(\gM_G^{\et})$ as an unramified $\GK$-representation via $\GKinfty/I_{K_\infty}\cong \GK/I_K$, we apply the full faithfulness of $\GKinfty$-restriction on crystalline representations \cite[Corollary~2.1.14]{kisin:fcrys} to conclude that $T^*_\Sig(\gM_G^{\et})$ is a maximal unramified $\GK$-quotient of $T_p(G)$. This shows that $T^*_\Sig(\gM_G^{\et})$ coincides with the quotient $T_p(G^{\et})$ of $T_p(G)$. Now Proposition~\ref{prop:Connd} follows.% from Tate's theorem \cite[(4.2)]{Tate:pDivGps} it is enough to show that
%
%One can show directly that $\gM\in\phimodBT\Sig$ is \'etale if and only if $T^*_\Sig(\gM)$ is unramified (using that $(\wh\Sig^{\ur})^{I_{K_\infty}} = W(\kbar)[[u]]$ and Lemma~\ref{lem:classifMult}). Now, note that a crystalline $\GK$-representation $V$ is unramified if and only if $V|_{\GKinfty}$ is unramified; indeed, via $\GKinfty/I_{K_\infty}\cong \GK/I_K$ one can extend  $V|_{\GKinfty}$ to an unramified $\GK$-representation (which is crystalline), and a crystalline $\GK$-representation is uniquely determined by its $\GKinfty$-restriction \cite[Corollary~2.1.14]{kisin:fcrys}. This proves the assertion on \'etale-ness in \eqref{prop:Connd:et}.
%
%The connected-ness assertion of  $p$-divisible groups now follows from Lemma~\ref{lem:ConndEtSeq}. By Cartier duality, we obtain the rest of the corollary for $p$-divisible groups. The assertions on finite flat group schemes follows by choosing an embedding into a $p$-divisible group.    
\end{proof}
%\begin{rmksub}
%Using Corollaries~\ref{cor:classifFF} and \ref{cor:FF} which will be proved later, one can adapt the above argument to prove the analogue of Proposition~\ref{prop:Connd} for finite flat group schemes of $p$-power order. 
%\end{rmksub}
%%
%%\begin{rmksub}\label{rmk:Connd}
%%Since the functor $\gM^*\circ T_p$ respects Cartier duality in the sense of \S\ref{par:DefEtPhiNilp}, we also have the Cartier dual statement of Proposition~\ref{prop:Connd}. In particular, $\gM^*\circ T_p$ takes the Cartier dual of the connected-\'etale sequence for $G^\vee$ to the exact sequence \eqref{eqn:ConndEtDual} for $\gM^*(G)$.
%%\end{rmksub}

The following proposition proves Proposition~\ref{prop:main} when $\gM$ is either $\vphi$-unipotent or of multiplicative type. %The proof in the  $\vphi$-unipotent case is not different from \cite[Proposition~1.2.7]{Kisin:2adicBT}, but let us give a proof for completeness.
\begin{prop}\label{prop:main2}
For any $\gM\in\phimodBT\Sig$, put $\M:=S\otimes_{\vphi,\Sig}\gM$. 
\begin{enumerate}
\item\label{prop:main2:Unip}
Assume  that $\gM$ is $\vphi$-unipotent, then the $\GKinfty$-map $\imath:T^*_\Sig(\gM) \hra T^*_{\st}(\M)$ defined in \eqref{eqn:TSigTqstLatt} is an isomorphism.
\item\label{prop:main2:Mult}
Assume that $\gM$ is of multiplicative type. If $p>2$ then the $\GKinfty$-map $\imath:T^*_\Sig(\gM) \hra T^*_{\st}(\M)$ defined in \eqref{eqn:TSigTqstLatt} is an isomorphism. If $p=2$, then the image of $\imath$ is precisely $2\tim T^*_{\st}(\M)$.
\end{enumerate}
\end{prop}
%%Recall that  $T^*_{\qst}(\M)$ and $T^*_{\st}(\M)$ define the same $\Zp$-lattice in $V:=V^*_{\cris}(D)$ by Lemma~\ref{lem:QstVsSt}. 
%Since $\imath$ respects embeddings into $V$ it follows that $T$ is $\GK$-stable sublattice in $T^*_{\st}(\M)$ if $\gM$ is either $\vphi$-unipotent or of multiplicative type, which proves Proposition~\ref{prop:main} in those cases.
\begin{proof}
The $\vphi$-unipotent case is exactly \cite[Proposition~1.2.7]{Kisin:2adicBT}, so it is left to handle the multiplicative-type case.
The injectivity of $\imath$ is clear. We will prove the proposition by looking at $(\imath \bmod p)$ and $(\imath \bmod4)$ when $p=2$.

For the proof, we may replace $K$ by $\wh K^{\ur}$ and $\gM$ by $W(\kbar)[[u]]\otimes_\Sig\gM$; indeed, by naturally identifying $V=V^*_{\cris}(D) \cong V^*_{\cris}(W(\kbar)\otimes_{W(k)}D)$ as $I_K$-representations, we do not change the relevant $\Zp$-lattices in $V$ and the map $\imath$. Now, by Lemma~\ref{lem:classifMult} it is enough to prove the proposition for $\gM = (\Sig,\PP(u)\vphi)$. The following claim concludes the proof of Proposition~\ref{prop:main2}\eqref{prop:main2:Mult}. 
\begin{claimsub}
Assume that $k=\kbar$ and let $\gM = (\Sig,\PP(u)\vphi)$. If $p>2$ then $(\imath\bmod p)$ is an isomorphism. If $p=2$ then $(\imath\bmod 2)$ is the zero map and $(\imath\bmod 4)$ is non-zero.
\end{claimsub}
Let $\alpha\in W(k)\starr$ be such that  $\alpha u^e \equiv \PP(u) \bmod{p}$, and choose $\beta\in W(k)\starr$ so that $\beta^{p-1} \equiv \alpha\iv \bmod p$. Then, for any $f\in\Hom_{\Sig,\vphi}(\gM/(p),\wh\Sig^{\ur}/(p))$ we have
\begin{equation*}
(f(\beta))^p= \beta^p(f(1))^p = u^e f(\beta).
\end{equation*}
If $f$ is non-zero then we have $f(\beta)^{p-1}= u^e$. By Lemma~\ref{lem:ker}, on the other hand, $\vphi:\wh\Sig^{\ur}/(p)\ra \Acris/(p)$ takes $f(\beta)$ to $0$ if and only if $p=2$. This proves the claim on $\imath\bmod p$.

Now, assume that $p=2$ and choose a $\Z_2$-basis $f\in\Hom_{\Sig,\vphi}(\gM,\wh\Sig^{\ur})$. Since  $f(\beta) \equiv  u^e \bmod{2}$, we can find $g\in\wh\Sig^{\ur}$ such that $f(\beta) = u^e + 2g$. Similarly, we can write  $\alpha\iv \PP(u) = u^e+2c(u)$ for some $c(u)\in\Sig$ of degree $<e$. Then we have
\begin{equation*}
\vphi(f(\beta)) \equiv \alpha\iv\PP(u)f(\beta) \equiv u^{2e}+ 2u^e(c(u) +  g) \bmod{4}.
\end{equation*}

Put $\wt f:= \imath(f) \in \Hom_{S,\vphi_1,\Fil^1}(\M,\Acris)$. Then by definition we have %$\wt f(\beta) = \vphi(f(\beta)) \in\Acris$, so
\begin{displaymath}
\wt f(\beta) \equiv u^{2e}+ 2u^e(c(u) +  g)  \equiv 2\big(\frac{u^{2e}}{2!} +u^ec(u)+u^e g\big) \bmod{4}.
\end{displaymath}
%This cannot be congruent to zero modulo $4$, since $u^{2e}/2!\notin \Sig^{\ur}/(2)$ but $u^ec(u)+u^e g\in \Sig^{\ur}/(2)$.
Since  $u^{2e}/2!+u^ec(u)+u^e g \not\equiv 0 \bmod 2$ (because $u^{2e}/2!\notin \Sig^{\ur}/(2)$ and $u^ec(u)+u^e g\in \Sig^{\ur}/(2)$), we have $\wt f (\beta)\not\equiv 0 \bmod 4$.
\end{proof}
\begin{lemsub}\label{lem:ker}
The kernel of the map $\vphi:\wh\Sig^{\ur}/(p) \ra \Acris/(p)$ is generated by $u^e$ where $e=\deg\PP(u)$.
\end{lemsub}
\begin{proof}
Recall that $R$ is isomorphic to the valuation ring of the completion of $k\llpar \nf\pi \rrpar^{\alg}$, and the subring $\wh\Sig^{\ur}/(p)\subset R $ is isomorphic to the valuation ring of $k\llpar u \rrpar^{\sep}$ (\emph{cf.} \S\ref{subsec:GKinftyRep}). Also one can directly check that $\Acris/(p)$ is the divided power envelop of $\rep$ with respect to the kernel of the $0$th projection $\bar\theta:\rep \thra \fo_{\Kbar}/(p)$.\footnote{This follows from the universal property of divided power envelop; \emph{cf.} \cite[Remark~3.20(8)]{Berthelot-Ogus}.} Since $\ker(\bar\theta)$ is principally generated by $\nf\pi^e$ \cite[Proposition~4.3.7]{wintenberger:NormFiels},  the kernel of the natural map $\rep\ra\Acris$ is principally generated by $\nf\pi^{pe}$. Now the lemma follows.
%By construction, the map $\vphi:\wh\Sig^{\ur}/(p) \ra \Acris/(p)$ is equal to the following  composite 
%\begin{equation*}
%\wh\Sig^{\ur}/(p) \xra{\vphi} \wh\Sig^{\ur}/(p) \hra R \ra \Acris/(p),
%\end{equation*}
%where the second and the third arrows are natural morphisms.
%
%Since any element $x\in \wh\Sig^{\ur}/(p) $ satisfies $x^N=\alpha u^n$ for some $\alpha\in (\wh\Sig^{\ur}/(p) )\starr$ and some integers $n$ and $N$, we have  $\vphi(x)=0$ in $\Acris$ if and only if $u^e\mid x$.
\end{proof}

The proposition below completes the proof of  Proposition~\ref{prop:main}, hence the proof of Theorem~\ref{thm:classif}:
\begin{prop}\label{prop:key}
For any $\gM\in\phimodBT\Sig$ and $\M:=S\otimes_{\vphi,\Sig}\gM$, the image of the map $\imath:T^*_\Sig(\gM)\ra T^*_{\st}(\M)$ is $\GK$-stable, where $\imath$ is defined in (\ref{eqn:TSigTqstLatt}). If $p>2$ then $\imath$ is an isomorphism, and if $p=2$ then the cokernel of $\imath$ is killed by $p=2$.
\end{prop}
Note that  a simpler proof when $p>2$ already appeared in the proof of  \cite[Theorem~2.2.7]{kisin:fcrys}.
\begin{proof}
From Proposition~\ref{prop:main2} and  the exactness results (Lemmas~\ref{lem:exactBreuil} and \ref{lem:TstExact}), it easily follows that $\imath$ is an isomorphism when $p>2$, and when $p=2$ we have 
\begin{equation*}
2T^*_{\st}(\M)\subseteq T^*_\Sig(\gM)\subseteq T^*_{\st}(\M).
\end{equation*}
%To show that $T^*_\Sig(\gM)$ is $\GK$-stable in $T^*_{\st}(\M)$, it suffices to show that the image of  $T^*_\Sig(\gM)$ in $T^*_{\st}(\M)/2T^*_{\st}(\M)$ is $\GK$-stable; i.e., for any $f\in T^*_\Sig(\gM)$ and $g\in \GK$, there exists $f'\in T^*_\Sig(\gM)$ such that 
%$
%g\tim\imath(f) \equiv \imath(f') \bmod{2T^*_{\st}(\M)} , 
%$
%where $g\tim\imath(f)$ is as defined in (\ref{eqn:TstGK}). Equivalently, we need to show that 
%$
%\big(g\tim\imath(f)\big)(1\otimes m) \in \vphi(\wh\Sig^{\ur}) + 2\Acris
%$ 
%for any $f\in T^*_\Sig(\gM)$, $m\in\gM$, and $g\in \GK$.
% Since $2\delta(1\otimes m)\in \vphi(\wh\Sig^{\ur})$ and $\wh\Sig^{\ur}$ is a saturated subring\footnote{$2W(\rep)\cap \vphi(\wh\Sig^{\ur}) = 2\vphi(\wh\Sig^{\ur})$} of $W(\rep)$, it suffices to show that 
% \begin{equation*}
% \big(g\tim \imath(f)\big)(1\otimes m)\in W(\rep)+2\Acris
% \end{equation*}
%  for any $f\in T^*_\Sig(\gM)$, $m\in\gM$, and $g\in \GK$.
To show that $T^*_\Sig(\gM)$ is $\GK$-stable in $T^*_{\st}(\M)$, it suffices to show that the image of  $T^*_\Sig(\gM)$ in $T^*_{\st}(\M)/2T^*_{\st}(\M)$ is $\GK$-stable; i.e., for any $f\in T^*_\Sig(\gM)$ and $g\in \GK$, there exists $f'\in T^*_\Sig(\gM)$ such that 
\begin{equation*}
g\tim\imath(f) \equiv \imath(f') \bmod{2T^*_{\st}(\M)} , 
\end{equation*}
where $g\tim\imath(f)$ is as defined in (\ref{eqn:TstGK}). 

Recall that we have the short exact sequence $0\ra\gM^u\ra\gM\ra\gM^m\ra0$ defined in \eqref{eqn:ConndEtDual}, which stays exact after applying $S\otimes_{\vphi,\Sig}(\cdot)$ by Lemma~\ref{lem:exactBreuil}.  We set $\M^u:=S\otimes_{\vphi,\Sig}\gM^u$ and $\M^m:=S\otimes_{\vphi,\Sig}\gM^m$. From Proposition~\ref{prop:main2}(\ref{prop:main2:Unip}), there exists $f'\in T^*_\Sig(\gM)$ such that  $\big(g\tim\imath(f)\big)|_{\M^u} = \imath(f')|_{\M^u}$. So the element $\delta:= g\tim\imath(f) - \imath(f')$ lies in $T^*_{\st}(\M^m)$. Now it suffices to show that $\delta\in 2T^*_{\st}(\M^m)$; i.e.,  $\delta(1\otimes y) \in 2\Acris$ for any $y\in\gM^m$. 

Pick $t\in \Acris$ as in Theorem~\ref{thm:Acris}. Since $\gM^m$ is multiplicative, we have 
\begin{equation*}
\delta(\M^m) \subset (W(\kbar)\tim S) \tim (t/2) \subset W(\rep)\tim (t/2).
\end{equation*} 
It follows  from Theorem~\ref{thm:Acris} that $W(\rep) \tim \frac{t}{2} \cap (W(\rep)+2\Acris) \subset 2\Acris$. (To see this,  assume that $a\tim\frac{t}{2}\in  W(\rep)+2\Acris$ for some $a\in W(\rep)$. Then one can find $b, c_n\in W(\rep)$ such that $a\tim\frac{t}{2} = b+2\sum_{n=1}^\infty c_n \gamma_n(\frac{t}{2})$ in $\Acris$. By Theorem~\ref{thm:Acris}, $b+ (-a+2c_1)\frac{T}{2} +2\sum_{n=2}^\infty c_n \gamma_n(\frac{T}{2}) \in \fa_{\cris}$ must be divisible by $[\nf\epsilon]-\sum_{n=0}^\infty \gamma_n(T)$. This forces that $a\in (2,[\nf\epsilon]-1)W(R)$, but since $[\nf\epsilon]-1 = -\sum_{i=1}^\infty 2^n \gamma_n(\frac{t}{2}) \in 2\Acris$ we conclude that $a\in 2\Acris$.) In particular, we are reduced to showing that $\delta(1\otimes y) \in W(\rep)+2\Acris$ for any $y\in\gM^m$. Now, recall that $\imath(f')(1\otimes x)\in \vphi(\wh\Sig^{\ur})$ for any $x\in\gM$ by construction, so it suffices to show that 
 \begin{equation}\label{eqn:keyClm}
 \big(g\tim \imath(f)\big)(1\otimes x)\in W(\rep)+2\Acris
 \end{equation}
  for any $f\in T^*_\Sig(\gM)$, $x\in\gM$, and $g\in \GK$.

With the notation as above, recall from (\ref{eqn:TstGK}) that
\begin{equation}\label{eqn:keyAcris}
\begin{aligned}
\big(g\tim \imath(f)\big)(1\otimes x) 
&= g  \Big(\sum_{i=0}^\infty \gamma_i(-t_g) \big(\imath(f)\circ N_{\M}^i\big)(1\otimes x)\Big)  \\
& = g \left( \imath(f)(1\otimes x)\right ) + g\Big(\sum_{i=1}^\infty \gamma_i(-t_g) \big(\imath(f)\circ N_{\M}^i\big) (1\otimes x)\Big),
\end{aligned}
\end{equation}
where $\gamma_i$ is the standard $i$th divided power. Inspecting the second row of (\ref{eqn:keyAcris}), the first term is in $W(\rep)$ by construction, and the second term is in $2\Acris$ as  $2$ divides $\gamma_i(t_g)$ for any $g\in\GK$ and $i>0$. This shows (\ref{eqn:keyClm}), hence the proposition.
\end{proof}

\subsection{Comparison with Dieudonn\'e crystals: proof of Proposition~\ref{prop:Dieudonne}(\ref{prop:Dieudonne:Comp})}\label{subsec:Dieudonne}\footnote{The author thanks Tong~Liu for providing his idea to improve the original argument.}
Let $G$ be a $p$-divisible group over $\fo_K$, $D$ the associated filtered $\vphi$-module, and  $\DD^*(G)$ the contravariantly associated Dieudonn\'e crystal. The $S$-module $\DD^*(G)(S)$ can be naturally viewed as a strongly divisible $S$-lattice in $S\otimes_{W(k)}D$. (See, for example, \cite[(A.1)--(A.2)]{kisin:fcrys} for more details.) One can actually recover the Dieudonn\'e crystal $\DD^*(G)$ from the strongly divisible $S$-module $\DD^*(G)(S)$.\footnote{This can be done by applying \cite[Theorem~6.6]{Berthelot-Ogus} to the unique differential operator given by Proposition~\ref{lem:N}.} %To simplify the notation, we let $\M_\DD:=\DD^*(G)(S)$ denote this strongly divisible $S$-lattice.

Let $\gM_G:=\gM^*(T_p(G))$ and we associated a strongly divisible $S$-lattice $S\otimes_{\vphi,\Sig}\gM_G$. The goal of this section is to prove the following lemma. 
\begin{lemsub}\label{lem:LattComp}
With the notation as above, we have  $\DD^*(G)(S)=S\otimes_{\vphi,\Sig}\gM_G$ as strongly divisible $S$-lattices in $S\otimes_{W(k)}D$. 
\end{lemsub}

The starting point of the proof is the following theorem of Faltings, which is the geometric ingredient of the proof.
\begin{thmsub}[{\cite[\S6, Theorem~7]{Faltings:IntegralCrysCohoVeryRamBase}}]\label{thm:FaltingsKisin}
There exists a natural $\GK$-equivariant \emph{injective} map $T_p(G) \hra T^*_{\st}(\DD^*(G)(S))$ which respects the embeddings into $T_p(G)[\ivtd p] \cong V^*_{\cris}(D)$. Furthermore, this map is an isomorphism if $p>2$. %or if $p=2$ and the Cartier dual $G^\vee$ is connected. 
%and the cokernel of the map is annihilated by $p$ if $p=2$.
\end{thmsub}
%\begin{proof}
%This is proved in Theorem~7 in \cite[\S6]{Faltings:IntegralCrysCohoVeryRamBase}, except that the map  $T_p(G) \hra T^*_{\st}(\DD^*(G)(S))$ is an isomorphism when $p=2$ and $G^\vee$ is connected. The latter is proved in \cite[Proposition~1.1.10]{Kisin:2adicBT} using Zink's theory of windows and display for connected $p$-divisible groups \cite{Zink:WindowsDisplayProgrMath195, Zink:DisplayFormalGpAsterisq278}.
%\end{proof}

To use Faltings's theorem to prove Lemma~\ref{lem:LattComp}, we need the following lemma:
\begin{lemsub}\label{lem:LattIncl}
Let $\M$ and $\M'$ be strongly divisible $S$-lattices in $S\otimes_{W(k)}D$, and assume that $\M\cong S\otimes_{\vphi,\Sig}\gM$ for some $\gM\in\phimodBT\Sig$. If $T^*_{\st}(\M) \subseteq T^*_{\st}(\M')$ as $\Zp$-lattices in $V:=V^*_{\cris}(D)$, then we have  $\M \supseteq \M'$.
\end{lemsub}
\begin{proof}
%Let us fix a $\Zp$-basis $\set{f_1,\cdots,f_d}$ of $T^*_{\st}(\M)$. %and view it as a $\Qp$-basis of $V_{\cris}(D)$ via the isomorphism. 
%Define $f:\M[\ivtd p]\ra (\Bcris)^d$ by $f(x) := ( f_1(x),\cdots, f_d(x))$ for any $x\in\M[\ivtd p]$. %Clearly, $f$ respects $\Fil^1$ and $\vphi_1$. 
%%
%From the comparison isomorphism $\Hom_{S}(\M,\Bcris)\liso \Bcris\otimes_{\Zp}T^*_{\st}(\M) $ (\emph{cf.} Lemma~\ref{lem:TstEmb}) it follows that $f$ induces a $\Bcris$-isomorphism $\Bcris\otimes_{S}\M \riso (\Bcris)^d$. In particular,   %$f:\M[\ivtd p]\ra (\Bcris)^d$ is injective.
%$f$ is injective.
Assume that $\M$ admits a $\Sig$-lattice $\gM$ as in the statement.
It suffices to show that if $x\in\M[\ivtd p]$ satisfies $f(x)\in \Acris$ for any $f\in T_{\st}^*(\M)$, then $m\in \M$. (The lemma follows by applying this to any $x\in\M'$.)

Assume to the contrary that  there exists $x\in\M[\ivtd p]\setminus \M$ such that $f(x)\in\Acris$ for any  $f\in T_{\st}^*(\M)$. Let $n$ be the smallest positive integer such that $p^nx\in\M$. Then we have $p^nx\in \M\setminus p\M$, and $f(p^nx) \in p^n\Acris$ for any $f\in T_{\st}^*(\M)$. But this contradicts to Lemma~\ref{lem:fbarinj} below.
\end{proof}
\begin{lemsub}\label{lem:fbarinj}
Suppose $\M = S\otimes_{\vphi,\Sig}\gM$ for some $\gM\in\phimodBT\Sig$, and let $y\in\M$. If $f(y)\in p\Acris$ for any $f\in T_{\st}^*(\M)$, then $y\in p\M$.
\end{lemsub}
\begin{proof}
By the ``multiplicative-unipotent sequence'' (\ref{eqn:ConndEtDual}), it suffices to handle the case when $\gM$ is either of multiplicative type or $\vphi$-unipotent. The case when $\gM$ is of multiplicative type is straightforward. (For example, one can use Lemma~\ref{lem:classifMult} to reduce to the rank-$1$ case.) When $\gM$ is $\vphi$-unipotent, the natural map $\iota:T^*_\Sig(\gM)\ra T^*_{\st}(\M)$ is an isomorphism, so it suffices to prove the same statement for $y\in\gM$ and $f\in T^*_\Sig(\gM)$.

Set $\ol\gM:=\gM/p\gM$ and $\ol\gN:=\bigcap_{\bar f\in T^*_\Sig(\ol\gM)}\ker(\bar f)\subseteq \ol\gM$. Since $\ol\gN$ is a saturated submodule, we have $\ol\gM/\ol\gN \in (\Mod/\Sig)^{\leqs1}$. (Note that $\Sig/(p) = k[[u]]$ is a discrete valuation ring.) By definition of $\ol\gN$ we have $T^*_\Sig(\ol\gM) = T^*_\Sig(\ol\gM/\ol\gN)$, which implies
\begin{equation*}\rank_{\Sig/(p)}\ol\gM  = \dim_{\Fp}T^*_\Sig(\ol\gM) =\dim_{\Fp}T^*_\Sig(\ol\gM/\ol\gN) = \rank_{\Sig/(p)}\ol\gM/\ol\gN.\end{equation*} Hence, $\ol\gN = 0$, which proves the lemma. 
\end{proof}

\begin{proof}[Proof of Lemma~\ref{lem:LattComp}]
Lemma~\ref{lem:LattComp} is easy for \'etale  and multiplicative-type $p$-divisible groups. (Indeed, it essentially boils down to handling $G= \Gm[p^\infty]$ or $\nf{\Qp/\Zp}$, in which case one can explicitly compute and compare both sides of the equality.) %Now by connected-\'etale sequence and Proposition~\ref{prop:Connd}, we only need to verify  Lemma~\ref{lem:LattComp} for the Cartier duals of connected $p$-divisible groups to handle the general case. 

To prove Lemma~\ref{lem:LattComp} it is enough to show an inclusion $S\otimes_{\vphi,\Sig}\gM_{G}\supset \DD^*(G)(S)$ for all $p$-divisible group $G$. Indeed, we can apply this inclusion to the Cartier dual $G^\vee$ and obtain the inclusion of the opposite direction. (Note that Cartier duality corresponds to $S$-linear duality on both sides by \S\ref{par:DefEtPhiNilp} and \cite[\S5.3]{Berthelot-Breen-Messing:DieudonneII}.)

To show this inclusion, it suffices to show $T^*_{\st}(S\otimes_{\vphi,\Sig}\gM_{G})\subset T^*_{\st}(\DD^*(G)(S))$ by Lemma~\ref{lem:LattIncl}. Let $G$ be any $p$-divisible group over $\fo_K$ and consider the multiplicative-unipotent sequence $(G^\bullet):\, 0\ra G^m \ra G \ra G^u\ra 0$; i.e., $G^m$ is of multiplicative type and $(G^u)^\vee$  is connected. %To simplify the notation, set $\gM_{G^\bullet}:=\gM^*(T_p(G^\bullet))$. 
By Proposition~\ref{prop:Connd}, the sequence $\gM_{G^\bullet}:=\gM^*(T_p(G^\bullet))$ is short exact and coincides with the ``multiplicative-unipotent sequence'' \eqref{eqn:ConndEtDual} for $\gM_G$. % $S\otimes_{\vphi,\Sig)(\gM^\bullet):\,  0 \ra \M^u \ra \M \ra \M^m \ra 0$ and $\DD^*(G^\bullet)(S):\, 0\ra \M^u_\DD \ra \M_\DD \ra \M^m_\DD \ra 0$. 

We already observed $T^*_{\st}(S\otimes_{\vphi,\Sig}\gM_{G^m}) = T^*_{\st}(\DD^*(G^m)(S))$.
Applying Proposition~\ref{prop:main2}(\ref{prop:main2:Unip}) and Theorem~\ref{thm:FaltingsKisin}, we obtain that $T^*_{\st}(S\otimes_{\vphi,\Sig}\gM_{G^u}) = T_p(G^u)$ is contained in $T^*_{\st}(\DD^*(G^u)(S))$ as  $\Zp$-lattices in $V_p(G^u)$. Using the exactness of $T^*_{\st}$ (Lemma~\ref{lem:TstExact}), one easily obtains the desired inclusion. 
\end{proof}

\subsection{Exactness}\label{subsec:exactness}
We now give a proof of Proposition~\ref{prop:Dieudonne}\eqref{prop:Dieudonne:Exact} and Corollary~\ref{cor:classifFF}, using compatibility with crystalline Dieudonn\'e theory (Proposition~\ref{prop:Dieudonne}\eqref{prop:Dieudonne:Comp}) proved in \S\ref{subsec:Dieudonne}. %Indeed, we also prove stronger exactness (Proposition~\ref{prop:exact}) that can be useful in practice.

The following lemma proves Proposition~\ref{prop:Dieudonne}\eqref{prop:Dieudonne:Exact} 
\begin{lemsub}\label{lem:exact}
Let $(G^\bullet):\ 0\ra G' \ra G \ra G'' \ra 0$ be a short exact sequence of $p$-divisible groups over $\fo_K$, and put $(\gM^\bullet):=\gM^*(G^\bullet)$.
\begin{enumerate}
\item
\label{lem:exact:GenFib}
The sequence $T_p(G^\bullet)$ of Tate modules is short exact if and only if $\fo_\Eps\otimes_\Sig(\gM^\bullet)$ is short exact. %, where $\fo_\Eps$ is the $p$-adic completion of $\Sig[\ivtd u]$.
\item
\label{lem:exact:exact}
If $(G^\bullet)$ is an exact sequence of $p$-divisible groups over $\fo_K$, then $(\gM^\bullet)$ is an exact sequence.
\end{enumerate}
\end{lemsub}
\begin{proof}
By definition we have natural $\GKinfty$-isomorphisms $T_p(G^\bullet)\cong T^*_\Sig(\gM^\bullet)$ and the latter is precisely the sequence of \'etale $\vphi$-modules associated to $T_p(G^\bullet)$ by \cite[\S{B}.1.8.3]{fontaine:grothfest}. Now the statement \eqref{lem:exact:GenFib} directly follows from the exactness property of \'etale $\vphi$-modules. See \cite[\S{A}.1.2]{fontaine:grothfest} for more details.

Assuming that $(G^\bullet)$ is an exact sequence of $p$-divisible groups over $\fo_K$, the left exactness of $(\gM^\bullet)$ follows from \eqref{lem:exact:exact} since the natural map $\gM \ra \fo_\Eps\otimes_\Sig\gM$ is injective for any $\gM\in\phimodBT\Sig$. By Nakayama lemma the right exactness could be checked mod $u$, and by the compatibility with crystalline Dieudonn\'e theory the sequence $(\gM^\bullet)\otimes_{\Sig,\vphi} \Sig/(u)$ is isomorphic to the sequence $D^*((G^\bullet)_k)$, where $D^*(G_k)$ denote the contravariant Dieudonn\'e module for a $p$-divisible group $G_k$ over $k$.\footnote{For any $p$-divisible group $G$ over $\fo_K$, we have $\DD^*(G)(S)\otimes_S W(k)\cong D^*(G_k)$, hence the desired isomorphism follows from Proposition~\ref{prop:Dieudonne}(\ref{prop:Dieudonne:Comp}).}  This shows the desired exactness.
\end{proof}

\parag
\label{par:FinFl}
Since we have proved Theorem~\ref{thm:classif} and Proposition~\ref{prop:Dieudonne}, Corollary ~\ref{cor:classifFF} also follows possibly except the exactness assertion (as remarked in the proof of Corollary ~\ref{cor:classifFF}). We now show a slightly stronger assertion than the exactness asserted in Corollary~\ref{cor:classifFF}, which can be useful in practice.
 
Consider a (not necessarily exact) sequence $(G^\bullet):\ 0\ra G' \ra G \ra G'' \ra 0$, 
%\begin{equation*}
%(G^\bullet):\ 0\ra G' \ra G \ra G'' \ra 0,
%\end{equation*}
where $G$, $G'$, and $G''$ are either $p$-power order finite flat group schemes or $p$-divisible groups over $\fo_K$. (From now on, all finite flat group schemes are assumed to be of $p$-power order.) Let $(\gM^\bullet):\ 0\ra\gM''\ra\gM\ra\gM'\ra0$ denote the (not necessarily exact) sequence of  $(\vphi,\Sig)$-modules corresponding to $(G^\bullet)$. If  $G'$ is a finite flat group scheme and $G$ is a $p$-divisible group then we define the map $\gM \ra \gM'$ as follows: we choose a sufficiently large $n$ so that $G[p^n]$ contains the image of $G'$, and consider $\gM \thra \gM/(p^n) \ra \gM'$ where the second map is induced from the map $G'\ra G[p^n]$ which factors $G'\ra G$. Clearly the map $\gM \ra \gM'$ is independent of the choice of $n$.
\begin{propsub}\label{prop:exact}
We use the same notation as above.
\begin{enumerate}
\item
\label{prop:exact:GenFib}
The sequence $(G^\bullet)_K$ of the generic fibers is exact if and only if $\fo_\Eps\otimes_\Sig(\gM^\bullet)$ is exact.
\item
\label{prop:exact:exact}
The sequence $(G^\bullet)$ is exact if and only if $(\gM^\bullet)$ is exact.
\end{enumerate}
\end{propsub}
\begin{proof}
The proof of Lemma~\ref{lem:exact}\eqref{lem:exact:GenFib} can easily be adapted to prove  \eqref{prop:exact:GenFib} using the natural $\GKinfty$-isomorphism $T^*_\Sig(\gM_H)\cong H(\Kbar)$ when $H$ is a finite flat group scheme. (When $G'$ is a finite flat group scheme and $G$ is a $p$-divisible group, we need the exact sequence \eqref{eqn:TSigRedn}.)

The ``only if'' direction of \eqref{prop:exact:exact} is already proved except when $G'$, $G$, and $G''$ are finite flat group schemes. This case can be handled via a similar argument to the proof of Lemma~\ref{lem:exact}\eqref{lem:exact:exact}; the left exactness can be deduced from \eqref{prop:exact:GenFib} and the injectivity of maps $(\gM^\bullet)\hra \fo_\Eps\otimes_\Sig(\gM^\bullet)$, and the right exactness can be checked by reducing $(\gM^\bullet)$ modulo $(u)$ and comparing it with the sequence of contravariant Dieudonn\'e modules for the special fibers $(G^\bullet)_k$.

The proof of the ``if'' direction can be reduced to the case when $G'$, $G$, and $G''$ are finite flat group schemes, which follows from the claim below:
\begin{nclma} %\label{clm:exact}
If $\mathfrak{i}:\gM\ra\gM'$ is surjective, then $i:G'\ra G$ is a closed immersion.
\end{nclma}
%\begin{nclmb} 
%If $\mathfrak{j}:\gM''\ra\gM$ is injective, then $j:G\ra G''$ is an epimorphism.
%\end{nclmb}
Granting this claim, we have $\gM_{G/G'} = \ker(\mathfrak i)=\gM''$ as submodules of $\gM$, so the natural map $G/G' \ra G''$ given by the universal property of quotient is an isomorphism. This proves the ``if'' direction of Proposition~\ref{prop:exact}(\ref{prop:exact:exact}).

Let us prove Claim~A. Let $H\subset G$ be the scheme-theoretic image of $i$, and $\gM_H$ the quotient of $\gM$ corresponding to $H$. Since $i$ factors through $H$, it follows that $\mathfrak{i}$ factors through $\gM_H$ and we obtain a natural surjective map $\gM_H\thra \gM'$. This map is an isomorphism if and only if  $i$ is a closed immersion. Since both the source and the target have no non-zero  $u$-torsion, it suffices to show that the map $\fo_\Eps\otimes\mathfrak{i}:\fo_\Eps\otimes_\Sig\gM_H\thra\fo_\Eps\otimes_\Sig \gM'$ an isomorphism. 
To show this, observe that the natural map $G'_K \ra H_K$ induced by $i$ on the generic fiber is an isomorphism by Proposition~\ref{prop:exact}\eqref{prop:exact:GenFib}. 
%Applying   Proposition~\ref{prop:exact}\eqref{prop:exact:GenFib} once again\footnote{One deduces from Proposition~\ref{prop:exact}\eqref{prop:exact:GenFib} that $\fo_\Eps\otimes\mathfrak{i}$ is a surjective map between modules with same length, hence it is an isomorphism. This also follows from the anti-equivalence of categories proved in \cite[\S{A}.1.2]{fontaine:grothfest}.}, 
By the theory of \'etale $\vphi$-modules \cite[\S{A}.1.2]{fontaine:grothfest}
we conclude that the natural surjective map $\fo_\Eps\otimes\mathfrak{i}$ is an isomorphism.\footnote{Alternatively, one deduces from Proposition~\ref{prop:exact}\eqref{prop:exact:GenFib} that $\fo_\Eps\otimes\mathfrak{i}$ is a surjective map between modules with same length, hence it is an isomorphism.} This concludes the proof of Claim~A (therefore, of Proposition~\ref{prop:exact}).
\end{proof}

\bibliography{bib}

\def\cprime{$'$}
\providecommand{\bysame}{\leavevmode\hbox to3em{\hrulefill}\thinspace}
\providecommand{\MR}{\relax\ifhmode\unskip\space\fi MR }
\providecommand{\MRhref}[2]{%
  \href{http://www.ams.org/mathscinet-getitem?mr=#1}{#2}
}
\providecommand{\href}[2]{#2}
\begin{thebibliography}{BBM82}

\bibitem[BBM82]{Berthelot-Breen-Messing:DieudonneII}
Pierre Berthelot, Lawrence Breen, and William Messing, \emph{Th{\'e}orie de
  {D}ieudonn{\'e} cristalline. {II}}, Lecture Notes in Mathematics, vol. 930,
  Springer-Verlag, Berlin, 1982. \MR{667344 (85k:14023)}

\bibitem[BK]{BlockKato:Dieudonne}
Spencer Bloch and Kazuya Kato, \emph{{$p$}-divisible groups and {Dieudonn\'e}
  crystals}, unpublished., 48 pages.

\bibitem[BO78]{Berthelot-Ogus}
Pierre Berthelot and Arthur Ogus, \emph{Notes on crystalline cohomology},
  Princeton University Press, Princeton, N.J., 1978. \MR{MR0491705 (58
  \#10908)}

\bibitem[Bre97]{Breuil:GriffithsTransv}
Christophe Breuil, \emph{Repr\'esentations {$p$}-adiques semi-stables et
  transversalit\'e de {G}riffiths}, Math. Ann. \textbf{307} (1997), no.~2,
  191--224. \MR{MR1428871 (98b:14016)}

\bibitem[Bre98]{breuil:GpSchNormField}
\bysame, \emph{Sch{\'e}mas en groupes et corps des normes}, Preprint,\hfill\\
  {\tt www.ihes.fr/\~{}breuil/PUBLICATIONS/groupesnormes.pdf} (1998).

\bibitem[Bre00]{Breuil:GrPDivGrFiniModFil}
\bysame, \emph{Groupes {$p$}-divisibles, groupes finis et modules filtr\'es},
  Ann. of Math. (2) \textbf{152} (2000), no.~2, 489--549. \MR{MR1804530
  (2001k:14087)}

\bibitem[Bre02]{Breuil:IntegralPAdicHodgeThy}
\bysame, \emph{Integral {$p$}-adic {H}odge theory}, Algebraic geometry 2000,
  {A}zumino ({H}otaka), Adv. Stud. Pure Math., vol.~36, Math. Soc. Japan,
  Tokyo, 2002, pp.~51--80. \MR{MR1971512 (2004e:11135)}

\bibitem[dJ95]{dejong:crysdieubyformalrigid}
A.~J. de~Jong, \emph{Crystalline {D}ieudonn{\'e} module theory via formal and
  rigid geometry}, Inst. Hautes {\'E}tudes Sci. Publ. Math. (1995), no.~82,
  5--96 (1996). \MR{1383213 (97f:14047)}

\bibitem[Fal99]{Faltings:IntegralCrysCohoVeryRamBase}
Gerd Faltings, \emph{Integral crystalline cohomology over very ramified
  valuation rings}, J. Amer. Math. Soc. \textbf{12} (1999), no.~1, 117--144.
  \MR{MR1618483 (99e:14022)}

\bibitem[Fon90]{fontaine:grothfest}
Jean-Marc Fontaine, \emph{Repr\'esentations {$p$}-adiques des corps locaux.
  {I}}, The Grothendieck Festschrift, Vol.\ II, Progr. Math., vol.~87,
  Birkh{\"a}user Boston, Boston, MA, 1990, pp.~249--309. \MR{MR1106901
  (92i:11125)}

\bibitem[Fon94]{fontaine:Asterisque223ExpII}
\bysame, \emph{Le corps des p\'eriodes {$p$}-adiques}, Ast\'erisque (1994),
  no.~223, 59--111, With an appendix by Pierre Colmez, P\'eriodes $p$-adiques
  (Bures-sur-Yvette, 1988). \MR{MR1293971 (95k:11086)}

\bibitem[Kis06]{kisin:fcrys}
Mark Kisin, \emph{Crystalline representations and {$F$}-crystals}, Algebraic
  geometry and number theory, Progr. Math., vol. 253, Birkh{\"a}user Boston,
  Boston, MA, 2006, pp.~459--496. \MR{MR2263197 (2007j:11163)}

\bibitem[Kis09a]{Kisin:2adicBT}
\bysame, \emph{Modularity of 2-adic {B}arsotti-{T}ate representations}, Invent.
  Math. \textbf{178} (2009), no.~3, 587--634. \MR{MR2551765}

\bibitem[Kis09b]{Kisin:ModuliGpSch}
\bysame, \emph{Moduli of finite flat group schemes and modularity}, Ann. of
  Math. (2) \textbf{170} (2009), no.~3, 1085--1180.

\bibitem[Kis10]{Kisin:IntModelAbType}
\bysame, \emph{Integral models for {S}himura varieties of abelian type}, J.
  Amer. Math. Soc. \textbf{23} (2010), no.~4, 967--1012.

\bibitem[Kud04]{Kudla:EisensteinOutline}
Stephen~S. Kudla, \emph{Special cycles and derivatives of {E}isenstein series},
  Heegner points and {R}ankin {$L$}-series, Math. Sci. Res. Inst. Publ.,
  vol.~49, Cambridge Univ. Press, Cambridge, 2004, pp.~243--270. \MR{2083214
  (2005g:11108)}

\bibitem[Lan76]{Langlands:Jugendtraum}
R.~P. Langlands, \emph{Some contemporary problems with origins in the
  {J}ugendtraum}, Mathematical developments arising from {H}ilbert problems
  ({P}roc. {S}ympos. {P}ure {M}ath., {V}ol. {XXVIII}, {N}orthern {I}llinois
  {U}niv., {D}e {K}alb, {I}ll., 1974), Amer. Math. Soc., Providence, R. I.,
  1976, pp.~401--418. \MR{0437500 (55 \#10426)}

\bibitem[Lau10a]{Lau:GalRep}
Eike Lau, \emph{Displayed equations for {Galois} representations}, Preprint,
  {\tt http://arxiv.org/abs/1012.4436v1} (2010).

\bibitem[Lau10b]{Lau:2010fk}
\bysame, \emph{A relation between dieudonne displays and crystalline dieudonne
  theory}, Preprint, {\tt http://arxiv.org/abs/1006.2720v1} (2010).

\bibitem[Liu08a]{Liu:LattFiltPhiNMod}
Tong Liu, \emph{Lattices in filtered {$(\varphi,N)$}-modules}, Preprint (2008).

\bibitem[Liu08b]{Liu:StronglyDivLattice}
\bysame, \emph{On lattices in semi-stable representations: a proof of a
  conjecture of {B}reuil}, Compos. Math. \textbf{144} (2008), no.~1, 61--88.
  \MR{MR2388556}

\bibitem[Liu10]{Liu:Classif}
\bysame, \emph{The correspondence between {B}arsotti-{T}ate groups and {K}isin
  modules when {$p=2$}}, Preprint (2010).

\bibitem[MM74]{Mazur-Messing}
B.~Mazur and William Messing, \emph{Universal extensions and one dimensional
  crystalline cohomology}, Lecture Notes in Mathematics, Vol. 370,
  Springer-Verlag, Berlin, 1974. \MR{0374150 (51 \#10350)}

\bibitem[MS11]{MadapusiPera:Thesis}
Keerthi Madapusi~Sampath, \emph{Toroidal compactifications of integral
  canonical models of {S}himura varieties of {H}odge type}, Ph.D. thesis,
  University of Chicago, 2011.

\bibitem[Ray74]{raynaud:GpSch}
Michel Raynaud, \emph{Sch\'emas en groupes de type {$(p,\dots, p)$}}, Bull.
  Soc. Math. France \textbf{102} (1974), 241--280. \MR{54 \#7488}

\bibitem[Tat67]{Tate:pDivGps}
John~T. Tate, \emph{{$p$-divisible} groups.}, Proc. Conf. Local Fields
  (Driebergen, 1966), Springer, Berlin, 1967, pp.~158--183. \MR{MR0231827 (38
  \#155)}

\bibitem[Vas07]{Vasiu:GoodRedn2}
Adrian Vasiu, \emph{Good reductions of shimura varieties of hodge type in
  arbitrary unramified mixed characteristic, part {II}}, Preprint {\tt
  http://arxiv.org/abs/0712.1572} (2007).

\bibitem[Win83]{wintenberger:NormFiels}
Jean-Pierre Wintenberger, \emph{Le corps des normes de certaines extensions
  infinies de corps locaux; applications}, Ann. Sci. \'Ecole Norm. Sup. (4)
  \textbf{16} (1983), no.~1, 59--89. \MR{MR719763 (85e:11098)}

\bibitem[Zin01]{Zink:DieudonnePDivGpCFT}
Thomas Zink, \emph{A {D}ieudonn{\'e} theory for {$p$}-divisible groups}, Class
  field theory---its centenary and prospect ({T}okyo, 1998), Adv. Stud. Pure
  Math., vol.~30, Math. Soc. Japan, Tokyo, 2001, pp.~139--160. \MR{MR1846456
  (2002h:14075)}

\end{thebibliography}
\bibliographystyle{amsalpha}

\end{document}